%% file: Witt_vectorsNEW.tex
\documentclass[10pt]{article}
\input{lang_en.tex}

\input{hat.tex}

\usepackage{tikz-cd}
\usepackage{array}
\usepackage{tikz}
\usetikzlibrary{decorations.pathmorphing}
\usepackage[alphabetic]{amsrefs}
\usepackage{hyperref}
\usepackage{bbm}
\usepackage{authblk}
\usetikzlibrary{babel}

\numberwithin{equation}{section}

\crefalias{prop}{proposition}

\newcommand{\F}{\mathrm{F}}

\newcommand{\LF}{\mathrm{LF}}
\newcommand{\CoAlg}{\mathrm{CoAlg}}

\newcommand{\CartAlg}{\mathrm{CartCAlg}}

\newcommand{\RL}{\mathrm{RL}}

\newcommand{\fg}{\mathrm{fg}}
\newcommand{\der}{\mathrm{der}}

\newcommand{\Spt}{\mathrm{Spt}}
\newcommand{\Spc}{\mathrm{Spc}}

\newcommand{\can}{\mathrm{can}}

\newcommand{\free}{\mathrm{free}}

\newcommand{\RMod}{\mathrm{RMod}}

\newcommand{\op}{\mathrm{op}}

\newcommand{\B}{\mathrm{B}}

\newcommand{\forget}{\mathrm{oubl}}

\newcommand{\LSym}{\mathrm{LSym}}

\newcommand{\ext}{\mathrm{ext}}
\newcommand{\D}{\mathrm{D}}
\newcommand{\apoly}{\mathrm{addpoly}}
\newcommand{\epoly}{\mathrm{excpoly}}
\newcommand{\poly}{\mathrm{poly}}

\newcommand{\Cart}{\mathrm{Cart}}
\newcommand{\DAlg}{\mathrm{DAlg}}
\newcommand{\DCartAlg}{\mathrm{CartDAlg}}
\newcommand{\CartCAlg}{\mathrm{CartCAlg}}

\newcommand{\W}{\mathrm{W}}

\newcommand{\LT}{\mathrm{LT}}

\newcommand{\C}{\mathrm{C}}

\newcommand{\T}{\mathrm{T}}

\newcommand{\Barr}{\mathrm{Bar}}

\DeclareMathOperator{\Id}{Id}

\mathchardef\mdef="2D 

\begin{document}

	\title{\textbf{Witt vectors and $\delta$-Cartier rings}}
	\date{}
	\author{Kirill Magidson}
	\affil{\small{Northwestern University}}
	\affil{\emph{\small{kirill.magidson@northwestern.edu}}}
	\maketitle
	
	\begin{abstract}
		We give a universal property of the construction of the ring of $p$-typical Witt vectors of a commutative ring, endowed with Witt vectors Frobenius and Verschiebung, and generalize this construction to the derived setting. We define an $\infty$-category of $p$-typical derived $\delta$-Cartier rings and show that the derived ring of $p$-typical Witt vectors of a derived ring is naturally an object in this $\infty$-category. Moreover, we show that for any prime $p$, the formation of the derived ring of $p$-typical Witt vectors gives an equivalence between the $\infty$-category of all derived rings and the full subcategory of all derived $p$-typical $\delta$-Cartier rings consisting of $V$-complete objects. 
	\end{abstract}

	\tableofcontents

	\section{Introduction.}

	\subsection{Motivating Witt vectors.}
	
	Let $A$ be a commutative ring, and fix a prime number $p$. The ring of  \textbf{($p$-typical) Witt vectors} $\W(A)$ is a fundamental algebraic construction which plays an important role in algebra, algebraic geometry and homotopy theory. Historically, the construction of the ring of Witt vectors was first introduced in the proof of the following classical theorem (see \cite[Chapter II, Theorem 5]{Ser}).
	
	\begin{Theor}\label{classical_theorem}
	Let $k$ be a perfect field of characteristic $p$. There is a unique $p$-adically complete ring $\W(k)$ with an isomorphism $F: \W(k) \xymatrix{\ar[r]^-{\sim}&}\W(k)$ such that $\W(k)/p \simeq k$, and the map $F$ reduces to the Frobenius map $\varphi: k \rightarrow k$.	\end{Theor}

	The construction of $\W(k)$ in this special case of a general functorial construction for any ring $A$.
	
	\begin{construction}\label{ghost_classical_def}
	For any commutative ring, let $\W(A):=\prod_{n\geq 0}A$ as a set, and let $\xymatrix{w_{n}: \prod_{n\geq 0} A \rightarrow A}$ be the map defined by the formula
	
	$$
	w_{n}(a_{0},a_{1},...):=a_{0}^{p^{n}}+pa_{1}^{p^{n-1}}+...+p^{n}a_{n}.
	$$
	There exists a unique commutative ring structure on $\W(A)$ such that the map $$\xymatrix{       \W(A) \ar[rr]^-{(w_{0},w_{1},...)}&&} \prod_{n\geq 0}A    $$
	is a natural transformation of functors on the category of rings (\cite[Chapter II, Theorem 7]{Ser}).
	\end{construction}

	\begin{construction}
	There exists a unique \textbf{Witt vectors Frobenius} map $F: \W(A) \rightarrow \W(A)$ which is a ring homomorphism determined uniquely by the condition that $F  w_{n}=w_{n+1}$.
	\end{construction}
	
	\begin{ex}
	Let $A$ be a ring of characteristic $p$. Then one can see that the map $F: \W(A) \rightarrow \W(A)$ defined by the formula $F(a_{0},a_{1},...) =(a_{0}^{p},a_{1}^{p},...)$ is the Witt vectors Frobenius.
	\end{ex}
	
	One can show that for a general ring $A$, the Witt vectors Frobenius $F: \W(A) \rightarrow \W(A)$ can be written as $F(a)=a^{p}+p\delta(a)$. It is natural to introduce the following definition.

	\begin{defn}
	A \textbf{$\delta$-ring} $(A,\delta)$ is a ring $A$ endowed with a set-theoretic operation $\delta:A \rightarrow A$ such that the map $\varphi: A \rightarrow A$ defined by the formula $\varphi(a)=a^{p}+p\delta(a)$ is a ring homomorphism.
	\end{defn}

	Let $\delta\text{-}\CAlg_{\mathbb{Z}} $ be the category of $\delta$-rings. The forgetful functor $\delta\text{-}\CAlg_{\mathbb{Z}} \rightarrow \CAlg_{\mathbb{Z}}$ admits a left adjoint (free $\delta$-ring) and a right adjoint (cofree $\delta$-ring). The following perspective on Witt vectors is due to A.Joyal \cite{J85}.
	
	\begin{Theor}\label{delta_Joyal}
	The functor $\W: \CAlg_{\mathbb{Z}} \rightarrow \delta\text{-}\CAlg_{\mathbb{Z}}$ is the cofree $\delta$-ring functor.
	\end{Theor}
The approach gives a universal property of the construction $A\longmapsto( \W(A), F)$. However, the cofree $\delta$-ring property of the ring of Witt vectors does not take into account one important piece of data. Namely, there is a map $V: \W(A) \rightarrow \W(A)$ which acts as the shift $V(a_{0},a_{1},...)=(0,a_{1},a_{2},...)$, whose image coincides with the kernel of the projection map $\W(A) \rightarrow A$. The map $V$ is known as the \textbf{Witt vectors Verschiebung}. An explicit computation shows that the map $V$ satisfies the following relations

$$
FV(x)=px
$$
for any $x \in \W(A)$ and 

\begin{equation}\label{V_proj}
V(F(x)y)=xV(y), \:\: V(xF(y))=V(x)y
\end{equation}
for any $x, y \in \W(A)$. The last two relations can be rephrased as saying that the map $V$ is an ideal $V: F_* \W(A) \rightarrow \W(A)$, and then the first relation says that the composition of this ideal with the ring homomorphism $F$ is the principal ideal generated by $p$. The map $V$ and its properties is not explicit in the definition of $\W(A)$ as a cofree $\delta$-ring. One of the main questions this text answers is the following.

\begin{question}
For any ring $A$, consider the ring of $p$-typical Witt vectors $\W(A)$ as a triple $(\W(A), F,V)$ where $F: \W(A)\rightarrow \W(A)$ and $V: \W(A) \rightarrow \W(A)$ are the Frobenius and Verschiebung. In what category does the triple $(\W(A),F,V)$ naturally lie, and what is the universal property of the construction $A \longmapsto (\W(A),F,V)$?
\end{question}

More generally, we would like to generalize the construction $A \longmapsto (\W(A),F,V)$ to non-discrete, i.e. derived rings. In the next subsection, we will summarize our approach to this question and formulate the main results of the paper.

\subsection{Formulation of the main results.}

The basic algebraic structures used in this paper are \textbf{derived rings}, and derived rings endowed with various additional structures. The theory of derived rings is a generalization of the theory of ordinary commutative rings and is is closely related to the theory of simplicial-cosimplicial commutative rings. The formal definition is as follows. Let $\mathbb{Z}$ be the ring of integers and $\Mod_{\mathbb{Z}}$ the $\infty$-category of $\mathbb{Z}$-module spectra. It contains the full subcategory $\Mod_{\mathbb{Z},\geq 0} \subset \Mod_{\mathbb{Z}}$ of connective objects, i.e. the objects whose homotopy groups are concentrated in non-negative degrees. Recall that for any $\infty$-category $\C$, the $\infty$-category $\End(\C)$ of endofunctors has a monoidal structure, and a \textbf{monad} is a unital algebra object in $\End(\C)$.

\begin{defn}
Let $\CAlg_{\mathbb{Z}}^{\poly}$ be the category of commutative polynomial rings on finitely many variables and all commutative ring maps between them. The $\infty$-category of \textbf{connective derived rings} $\DAlg_{\mathbb{Z},\geq 0}$ is the sifted completion of $\CAlg_{\mathbb{Z}}^{\poly}$:

$$
\DAlg_{\mathbb{Z},\geq 0}:=\mathcal{P}_{\Sigma}(\CAlg_{\mathbb{Z}}^{\poly} ).
$$
The forgetful functor $\DAlg_{\mathbb{Z},\geq 0} \rightarrow \Mod_{\mathbb{Z},\geq 0}$ admits a left adjoint. Let $\LSym_{\mathbb{Z}}: \Mod_{\mathbb{Z},\geq 0}\rightarrow \Mod_{\mathbb{Z},\geq 0}$ be the resulting monad. It preserves sifted colimits, and in addition to that, it preserves certain totalisations. Using the latter property, one shows following the works \cite{BM19}, \cite{R}, \cite{BCN21} (see also the exposition in \cite{pdI}) that the monad $\LSym_{\mathbb{Z}}: \Mod_{\mathbb{Z},\geq 0} \rightarrow \Mod_{\mathbb{Z},\geq 0}$ uniquely extends to a sifted colimit preserving monad $\LSym_{\mathbb{Z}}: \Mod_{\mathbb{Z}} \rightarrow \Mod_{\mathbb{Z}}$. The $\infty$-category of \textbf{derived rings} is defined as the $\infty$-category of $\LSym_{\mathbb{Z}}$-algebras in $\mathbb{Z}$-module spectra

$$
\DAlg_{\mathbb{Z}}:=\Alg_{\LSym_{\mathbb{Z}}}(\Mod_{\mathbb{Z}}).
$$
\end{defn}

For any derived ring $A$, we will define a new derived ring $\W(A)$ by a universal property. Let $\overline{\varphi}\text{-}\DAlg_{\mathbb{Z}}$ be the $\infty$-category of derived rings $A$ endowed with an operation $\overline{\varphi}: A \rightarrow A/p$. There is a functor $\CAlg_{\mathbb{Z},\heartsuit} \rightarrow \overline{\varphi}\text{-}\CAlg_{\mathbb{Z},\heartsuit}$ endowing any commutative ring $A$ with the $p$-th power operation $\overline{\varphi}: A \rightarrow A/p, x \longmapsto \overline{x}^{p}$. Following A.Holeman \cite{H22}, there exists an extension of this functor to a functor on the $\infty$-category of derived rings $\Frob: \DAlg_{\mathbb{Z}} \rightarrow \overline{\varphi}\text{-}\DAlg_{\mathbb{Z}}$. For a derived ring $A$, we call the map $\overline{\varphi}: A \rightarrow A/p$ the \textbf{derived Frobenius}.

\begin{defn}
A \textbf{derived $\delta$-ring} $(A,\varphi)$ is a derived ring $A$ endowed with an endomorphism $\varphi: A \rightarrow A$ and a homotopy between the composite $\xymatrix{A \ar[r]^-{\varphi}& A \ar[r]^-{\can} & A/p }$ and the derived Frobenius. 
\end{defn}

\begin{ex}
The free $\delta$-ring on one generator $\mathbb{Z}\{x\}$ is free on infinitely many generators $\mathbb{Z}\{x\}\simeq \mathbb{Z}[x_{1},x_{2},...]$, and the Frobenius lift is given by the formula

$$
F(x_{i})=x_{i}^{p}+px_{i+1}.
$$
\end{ex}

Let $\delta\text{-}\DAlg_{\mathbb{Z}}$ be the $\infty$-category of derived $\delta$-rings. The forgetful functor $\delta\text{-}\DAlg_{\mathbb{Z}} \rightarrow \DAlg_{\mathbb{Z}}$ admits a left and right adjoint. We will show that the right adjoint \textbf{derived Witt vectors ring} functor $\W: \DAlg_{\mathbb{Z}} \rightarrow \delta\text{-}\DAlg_{\mathbb{Z}}$ coincides with the classical functor of Witt vectors on discrete rings, and is right-left extended from the full subcategory of discrete rings to all derived rings. Moreover, the derived Witt vectors ring $\W(A)$ carries an additional Verschiebung operation $V: \W(A) \rightarrow \W(A)$ which is subject to certain conditions listed in the definition below. The triple $(\W(A),F,V)$ is an example of a derived $\delta$-Cartier ring, a notion we will introduce below in two steps.

First, we will define the notion of a derived Cartier ring—an algebraic gadget which packages the maps $F,V$ on $\W(A)$, forgetting that it is a derived $\delta$-ring.

\begin{defn}
A \textbf{derived Cartier ring} is a triple $(A,F,V)$ consisting of the following data:

\begin{itemize}

\item A derived ring $A$;

\item A derived ring endomorphism $F: A \rightarrow A$.

\item An $A$-module map $V: F_{*}A \rightarrow A$ which is a derived ideal in the sense that the quotient $A/V:=\cofib(V)$ is a derived ring, and the natural quotient map $A \rightarrow A/V$ is a derived ring map;

\item A homotopy $$F\circ V \sim p $$ of derived ideals.

\end{itemize}
\end{defn}

\begin{rem}
Assume $A$ is discrete. The map $V: F_{*}A \rightarrow A$ is just an endomorphism $V: A \rightarrow A$ satisfying the relations \ref{V_proj}.
\end{rem}

Below we will give a description of the free discrete Cartier ring. 

\begin{construction}\label{add_V_intro}
Suppose $A$ is ring with an endomorphism $F: A \rightarrow A$, Let $A[V]:=\oplus_{n \geq 0} AV^{n}$ be the free $A$-module on infinitely many generators, whose elements we write in the form of finite sums $\sum_{n\geq 0} V^{n} a_{n}$. The$A$-module structure on $A[V]$ is determined by the formula

$$
a.V^{n}(b)= V^{n}(F^{n}(a) b)   ,
$$
for any $a, b\in A$. There is also an $A$-linear multiplication given by

$$
V^{n}(a)V^{m}(b) =  p^{m}V^{m}(F^{k}(a)b)
$$
for $n=m+k$, $k\geq 0$. Then $A[V]$ becomes a commutative associative $A$-algebra.

Define a linear map $F: A[V] \rightarrow A[V]$ by the formula

$$
F(\sum_{n\geq 0}V^{n}a_{n}) = F(a_{0})+p\bigl( \sum_{n\geq 0} V^{n}a_{n+1}\bigr).
$$
Then $F$ is a ring homomorphism, and the shift map $V$ is an ideal $V: F_{*}A \rightarrow A$ such that $FV=p$. Hence, $A[V]$  is a Cartier ring. It has the universal property that for any other Cariter ring $B$ and a ring homomorphism $f: A \rightarrow B$ intertwining the $F$ maps, there is a unique extension of $f$ to a Cartier ring map $A[V] \rightarrow B$.
\end{construction}

The notion of a derived $\delta$-Cartier ring is a refinement of the notion of a derived Cartier ring, where the endomorphism $F$ is now a derived $\delta$-structure, and a certain compatibility between $\delta$ and $V$ is satisfied. To explain this notion, let us apply Construction \ref{add_V_intro} to the case $A=\mathbb{Z}\{x\}$, the free $\delta$-ring on one generator.

\begin{ex}

Consider the Cartier ring $\mathbb{Z}\{x\}[V]$. The Cartier ring Frobenius $F: \mathbb{Z}\{x\}[V] \rightarrow \mathbb{Z}\{x\}[V]$. To see this, we consider the equatiion

$$
FV(a)=V(a)^{p}+p\delta(a)=pa,
$$
in which we have $V(a)^{p}=p^{p-1}V(a^{p})$, and hence there is a solution

\begin{equation}\label{delta_V_for}
\delta(a)=a-p^{p-2}V(a^{p}),
\end{equation}
  unique by $p$-torsion freeness of $\mathbb{Z}\{x\}[V]$.
\end{ex}

Essentially, the notion of a derived $\delta$-Cartier ring is designed in such a way that $\mathbb{Z}\{x\}[V]$ is the free such an object on one generator. We will now give a rigorous definition.
\begin{defn}
Let $\delta\text{-}\CartCAlg^{\poly} $ be the subcategory of all Cartier rings containing $\mathbb{Z}\{x\}[V]$ and closed under finite coproducts. The sifted completion $ \mathcal{P}_{\Sigma}(  \delta\text{-}\CartCAlg^{\poly}  )$ is monadic over $\Mod_{\mathbb{Z},\geq 0}$, and the monad has a unique sifted colimit preserving extension to a monad on $\Mod_{\mathbb{Z}}$. Denote this monad $\LSym^{\delta}[V]: \Mod_{\mathbb{Z}} \rightarrow \Mod_{\mathbb{Z}}$. We define the $\infty$-category of \textbf{derived $\delta$-Cartier rings} as the $\infty$-category of $\LSym^{\delta}[V]$-algebras:

$$
\delta\text{-}\Cart\DAlg:=\Alg_{\LSym^{\delta}[V]}(\Mod_{\mathbb{Z}}).
$$
\end{defn}

\begin{rem}
By construction, since the relation \ref{delta_V_for} holds for free $\delta$-Cartier rings, it holds for all $\delta$-Cartier rings (up to coherent homotopies for non-discrete ones). Let us explain the meaning of this relation. Suppose we have a general Cartier ring $(A,F,V)$ such that $F$ is a lift of Frobenius. Then we have $$F(Va)=(Va)^{p}+p\delta(Va)=pa .$$
Rewriting $(Va)^{p}=p^{p-1}V(a^{p})$, we get $$p\delta(Va)=pa-p^{p-1}V(a^{p}) .$$
If $A$ is $p$-torsion free, then we can divide by $p$ and obtain the relation \ref{delta_V_for} for $\delta(Va)$. Therefore, any $p$-torsion free Cartier ring for which $F$ is a lift of Frobenius, gives an example of a $\delta$-Cartier ring. But if $A$ has $p$-torsion, then the relation \ref{delta_V_for} holds only after multiplication with $p$, and hence it is not in general a $\delta$-Cartier ring. Hence the relation \ref{delta_V_for} distinguishes $\delta$-Cartier rings from general Cartier rings, for which Frobenius comes from a $\delta$-structure.
\end{rem}

The significance of the $\infty$-category of derived $\delta$-Cartier rings is due the observation that if $A$ is a discrete ring, then the ring of Witt vectors $\W(A)$ is a $\delta$-Cartier ring. Using right-left extension, the construction extends to a functor $\W: \DAlg_{\mathbb{Z}} \rightarrow \delta\text{-}\Cart\DAlg$.

\begin{Theor}\label{main_theorem}

Let $\delta\text{-}\Cart\DAlg$ be the $\infty$-category of derived $\delta$-Cartier rings, and $(-)/V: \delta\text{-}\Cart\DAlg \rightarrow \DAlg_{\mathbb{Z}}$ the functor which sends a derived $\delta$-Cartier ring $A$ to the quotient $A/V$ considered as a mere derived ring.

\begin{enumerate}

\item For every prime $p$, there exists an adjunction

$$
\begin{tikzcd}
\delta\text{-}\Cart\DAlg \arrow[rr, shift right = -.75ex, "(-)/V"] && \arrow[ll, shift right = -.75ex, "\W"]    \DAlg_{\mathbb{Z}}, \\
\end{tikzcd}
$$

\item The composition

$$
\xymatrix{   \DAlg_{\mathbb{Z}} \ar[r]^-{\W}& \delta\text{-}\Cart\DAlg \ar[r]^-{\forget} &\delta\text{-}\DAlg_{\mathbb{Z}}   }
$$
is the cofree derived $\delta$-ring functor.

\item

The adjunction restricts to an equivalence\footnote{An interesting feature of this equivalence is that the left-hand side depends on $p$, while the right-hand side does not. However, the functor $\W$ involved in the construction depends on $p$ because it depends on the derived Frobenius operation $A\rightarrow A/p$. In other words, we get different equivalences for different primes $p$.}
$$
\begin{tikzcd}
(-)/V: \delta \text{-}\widehat{\Cart}\DAlg \arrow[rr, shift right = -.75ex, swap, "\sim"] && \arrow[ll, shift right = -.75ex] \DAlg_{\mathbb{Z}}: \W.
\end{tikzcd}
$$
\end{enumerate}
\end{Theor}

Theorem \ref{main_theorem} is the main result of this text. It gives a universal property of the derived ring of Witt vectors, in which all structures known to be classically present on it, are packaged into a compound algebraic gadget.

As an application of this result, we show how the equivalence between derived $V$-complete $\delta$-Cartier rings and derived rings recovers the following theorem of B.Antieau \cite{A23}.
\begin{Theor}
Let $\DAlg_{\mathbb{F}_{p}}^{\perf}$ be the $\infty$-category of derived $\mathbb{F}_{p}$-algebras having the property that the derived Frobenius endomorphism $\varphi: A \rightarrow A$ is an equivalence, and $\delta\text{-}\DAlg_{\mathbb{Z}_{p}}^{\perf, \wedge} $ be the $\infty$-category of $p$-complete derived $\delta$-rings whose Frobenius lift is an equivalence. Then we have an equivalence

$$
\begin{tikzcd} (-)/p: \delta\text{-}\DAlg_{\mathbb{Z}_{p}}^{\perf, \wedge}   \arrow[rr, shift right = -.75ex, swap, "\sim"] && \arrow[ll, shift right = -.75ex]  \DAlg_{\mathbb{F}_{p}}^{\perf}: \W    \end{tikzcd}
$$
\end{Theor}
Note that the classical Theorem \ref{classical_theorem} which initiated the study of Witt vectors is a special case of this general theorem.

\section*{Notation and terminology.}

We use the language of $\infty$-categories and higher algebra developed by J.Lurie in \cite{HTT} and \cite{HA}. We freely use the notions of an \textbf{$\infty$-category}, a \textbf{functor between $\infty$-categories}, an \textbf{$\infty$-categorical adjunction} etc. The amount of higher categorical machinery used in the text, is limited to a very basic level. Essentially the only technical statements we use is the \textbf{adjoint functor theorem}, see \cite[Corollary 5.5.2.9 ]{HTT}, and the \textbf{Barr-Beck-Lurie monadicity theorem} formulated \cite[Theorem 4.7.3.5]{HA}. In general, most of our constructions are formulated in a homotopy invariant way. In particular, all limits and colimits are in general assumed to be in the homotopy sense. For better exposition, we often restrict our attention to some discrete computations, in which case we specify that we are working with discrete cochain complexes, algebras etc. 

We let $\Spc$ be the $\infty$-category of \textbf{spaces}, or \textbf{$\infty$-groupoids}, and $\Spt$ the $\infty$-category of \textbf{spectra}. For a ring spectrum $R$, we denote $\Mod_{A}$ the $\infty$-category of \textbf{$A$-module spectra}. For example, our notation for the derived $\infty$-category of abelian groups is $\Mod_{\mathbb{Z}}$, and we refer to its objects simply as "$\mathbb{Z}$-modules". In a stable $\infty$-category $\fcat C$, we denote $[1]$ and $[-1]$ suspension and loop functors respectively. We use homological notation for the $t$-structures. For example, $\Mod_{\mathbb{Z},\geq 0}$ is the $\infty$-category of connective $\mathbb{Z}$-modules, and $\Mod_{\mathbb{Z},\heartsuit}$ for the abelian category of abelian groups.

We often encounter pull-back diagrams of $\infty$-categories. Unless specified otherwise, these are taken in the $\infty$-category of large $\infty$-categories. 

For a pair of $\infty$-categories $\fcat C$ and $\fcat D$, we denote $\Fun(\fcat C, \fcat D)$ the $\infty$-category of functors from $\fcat C$ to $\fcat D$.

\section*{Acknowledgements.} I am thankful to Ben Antieau, Lukas Brantner, Chris Brav, Adam Holeman, Zhouhang Mao, Akhil Mathew, Thomas Nikolaus and Allen Yuan for various conversations related to the topics of the text. Akhil Mathew communicated to me that the notion of $\delta$-Cartier rings was recently independently studied by V.Drinfeld.

	\section{Preliminaries.}

	\subsection{Recollection on Witt vectors.}

	In this subsection we will review the classical theory of the ring of $p$-typical Witt vectors. This material is classical, and more details can be found, for instance, in \cite[Chapter IX, §1]{Bour} and \cite[Chapter 1]{Ill79}.
	
	\begin{construction}\label{Witt_ghost}
	Let $A$ be a commutative ring and fix a prime $p$ and $n \geq 0$. The set of \textbf{$p$-typical Witt vectors of length $n+1$}\footnote{We use the convention in which $\W_{1}(A)=A$.} is the product $$\W_{n+1}(A):= \prod_{i=0}^{n}A.$$ We write elements $a \in \W_{n+1}(A)$ as $n+1$-tuples $a=(a_{0},a_{1},...,a_{n})$ and refer to $a_{i}\in A$ as "coordinates" of the "Witt vector" $a$. To endow the set of Witt vectors of $A$ with a ring structure, one uses the \textbf{ghost components} map $$\W_{n+1}(A) \rightarrow \prod_{i=0}^{n}A$$ sending a tuple $(a_{0},a_{1},...,a_{n})$ to the tuple $(w_{0},w_{1},...,w_{n})$, where each $w_{k}$ is given by the formula 
 
 \begin{equation}\label{ghost_components}
 w_{k}= \sum_{i=0}^{k} p^{i}a_{i}^{p^{k-i}}.
 \end{equation}
The classical computation in the theory of Witt vectors says that for any $A$, there exists a unique ring structure on $\W_{n}(A)$, such that the maps $\{w_{i}\}_{i=0}^{n}: \W_{n+1}(A) \rightarrow \prod_{i= 0}^{n}A$ assemble into a natural transformation of endofunctors of the category of rings.
 
 \end{construction}
 
 \begin{rem}
 
There exist universal polynomials $S_{1}(a,b),S_{2}(a,b),...,S_{n}(a,b)$ and $P_{1}(a,b),P_{2}(a,b),...,P_{n}(a,b)$ governing the addition and multiplication laws $\W_{n}(A)$:

$$
a+b=(S_{1}(a,b),... \:,S_{n}(a,b) ),
$$

$$
ab=(P_{1}(a,b),... \:,P_{n}(a,b)).
$$
One can write down explicit formulas for these polynomials in low degrees, and in principle compute them inductively. For instance, $$S_{0}(a,b)=a_{0}+b_{0}, \:\:S_{1}(a,b)=a_{1}+b_{1}+\frac{1}{p}\bigl(a_{0}^{p}+b_{0}^{p}-(a_{0}+b_{0})^{p}\bigr).$$
and

$$
P_{0}(a,b)=a_{0}b_{0}, \:\:P_{1}(a,b)=a_{0}^{p}b_{1}+a^{p}_{1}b_{0}+pa_{1}b_{1}.
$$
 We denote $[-]: A \rightarrow \W_{n}(A)$ the Teichmüller map sending any $a\in A$ to the Witt vector $(a,0,...0)$. The Teichmüller map is unital and multiplicative, but not additive since $S_{1}([a],[b])\neq 0$, and thus $[a+b]\neq [a]+[b]$.
 \end{rem}

For any $n$, there is a \textbf{restriction map} $R_{n}:\W_{n+1}(A) \rightarrow \W_{n}(A)$, which simply forgets the last component in the direct product. This map is a ring homomorphism for any $n$.
 
 \begin{defn}\label{Witt_vectors_ring}
 The \textbf{ring of Witt vectors} of a commutative ring $A$ is the limit of finite length Witt vectors along restriction maps:
 
 $$
 \W(A) :=\underset{\leftarrow}{\lim} \:\W_{n}(A).
 $$
 \end{defn}
 
 \begin{ex}
	The most classical example of a computation of the ring of Witt vectors is $\W(\mathbb{F}_{p}) \simeq \mathbb{Z}_{p}$. Indeed, let $[-]:\mathbb{F}_{p} \rightarrow \mathbb{Z}_{p}$ be the Teichmüller representative map sending $a\in \mathbb{F}_{p}$ to the $p$-adic integer $(a,0,0,...) \in \mathbb{Z}_{p}$. Then we have a map $\W(\mathbb{F}_{p}) \rightarrow \mathbb{Z}_{p}$ sending a Witt vector $(a_{1},a_{2},...)$ to the series $\sum_{n\geq 0} [a_{n}]p^{n}$, and one can check that this map is a ring isomorphism.
	\end{ex}

 \begin{construction}
 
 Fix $ n \geq 1$. There exists a multiplicative \textbf{Witt vectors Frobenius} map $F_{n}: \W_{n+1}(A) \rightarrow \W_{n}(A)$. In terms of the ghost component maps, the Frobenius $F_{n}: \W_{n+1}(A) \rightarrow \W_{n}(A)$ is determined by the requirement that $$F_{n} w_{n}=w_{n+1}.$$ More concretely, the Frobenius map is given by the formula $F(x_{0},x_{1},...)=(f_{0},f_{1},...)$, where $f_{0}=a_{0}^{p}+pa_{1}$ and each $f_{i}$ is determined by iterated equations (see \cite[Equation 1.3.4]{Ill79})
 
 \begin{equation}\label{iterated_equations}
 f_{0}^{p^{i}}+pf_{1}^{p^{i-1}}+...+p^{i}f_{i}=a_{0}^{p^{i+1}}+pa_{1}^{p^{i}}+...+p^{i+1}a_{i+1}.
 \end{equation}
 
 \end{construction}
 
 \begin{rem}
 Notice that it follows from the equations \ref{iterated_equations} that $f_{i}(a)=a_{i}^{p}$ $\mod p$.  In particular, if $p=0$ in $A$, then the Frobenius $F: \W_{n+1}(A) \rightarrow \W_{n}(A)$ is given by the formula $$F(a_{0},a_{1},...,a_{n})=(a_{0}^{p},a_{2}^{p},...,a_{n-1}^{p}).$$

 \end{rem}

 \begin{construction}
 The \textbf{Verschiebung map} is an additive map $V_{n}: \W_{n}(A) \rightarrow \W_{n+1}(A)$ is defined as a "shift" $$V(a_{0},a_{1},...,a_{n}) = (0,a_{0},a_{1},...,a_{n}).$$
The following relation is satisfied as endomorphisms on $\W_{n}(A)$ (see \cite[Equation 1.3.7]{Ill79}):

$$
F_{n}V_{n}=p.
$$
	\end{construction}

	\begin{rem}
	
	In terms of the Verschiebung operation, any Witt vector $(a_{0},a_{1},...) \in \W(A)$ can be written as the infinite sum $\sum_{n\geq 0} V^{n}[a_{n}]$. Similarly, any element of the ring of finite Witt vectors of length $n$, $(a_{0},...,a_{n-1}) \in \W_{n}(A)$ can be written as a finite sum $\sum_{i=0}^{n-1} V^{i}([a_{i}])$.
	\end{rem}

	For every $n$, the operations $F_{n}:\W_{n+1}(A) \rightarrow \W_{n}(A)$ and $V_{n}: \W_{n}(A) \rightarrow \W_{n+1}(A)$ intertwine the restriction maps, and passing to the limit, we obtain a pair of endomorphisms $F: \W(A) \rightarrow \W(A)$ and $V: \W(A) \rightarrow \W(A)$ satisfying the relation $$FV=p.$$  The map $F$ is a ring homomorphism, and $V$ is merely an additive map.

It turns out that in general, the Witt vectors Frobenius $F: \W(A) \rightarrow \W(A)$ can be written in the form $F(x)=x^{p}+p\delta(x)$ for a unique map of sets $\delta: \W(A) \rightarrow \W(A)$. The most efficient way to see the existence of a canonical $\delta$-operation is by realizing that $\W(A)$ has a universal property of a cofree $\delta$-ring. This observation is due to A.Joyal \cite{J85}. We will now briefly review the theory of $\delta$-rings following the language and notation of the text \cite{BS22} of B.Bhatt and P.Scholze.

\begin{defn}
A \textbf{$\delta$-ring} is a pair $(A,\delta)$ where $A$ is an ordinary commutative ring, and $\delta: A \rightarrow A$ is an operation (a map of sets) satisfying the following equations:

$$
\delta(0)=\delta(1)=0\\
$$

$$
\delta(xy)=x^{p} \delta(y) + \delta(x)y^{p} + p\delta(x)\delta(y).\\
$$

$$
\delta(x+y) = \delta(x) + \delta(y) + \frac{x^{p}+y^{p}-(x+y)^{p}}{p}.
$$
A map of $\delta$-rings $(A,\delta) \rightarrow (B,\delta')$ is a ring map which commutes with $\delta$-operations. We denote $\delta\text{-}\CAlg_{\mathbb{Z},\heartsuit}$ the category of $\delta$-rings.
\end{defn}

\begin{rem}
Let $A$ be a $p$-torsion free ring. Then the data of a $\delta$-structure on $A$ is the same as the data of specifying a ring homomorphism $F: A \rightarrow A$ which lifts the Frobenius. For any such map, we can write $F(x)=x^{p}+p\delta(x)$ for a unique $\delta(x)$, and the map $x \longmapsto \delta(x)$ defines a $\delta$-operation. More generally, it is an observation of Bhatt-Scholze that even if the ring is not $p$-torsion free, the $\delta$-structure is still equivalent to a lift of Frobenius understood in the \emph{derived sense} (see \cite[Remark 2.5]{BS22} for more details). This is the perspective that we will adopt for dealing with derived $\delta$-rings.
\end{rem}

Let us give some examples of $\delta$-rings. 

\begin{ex}
The ring $\mathbb{Z}$ has a $\delta$-structure with the lift of Frobenius $F=\Id$.
\end{ex}

\begin{ex}\label{Z{x}}
The free $\delta$-ring on one generator $\mathbb{Z}\{x\}$ is the polynomial ring on infinitely many variables $\mathbb{Z}\{x\}\simeq \mathbb{Z}[x_{1}, x_{2},x_{3},...] $ with the lift of Frobenius defined on generators by the formulas (\cite[Lemma 2.11]{BS22}) \begin{equation}\label{F_on_Z{x}} F(x_{n})=x^{p}_{n}+px_{n+1}.\end{equation}
\end{ex}

\begin{ex}
The polynomial ring $\mathbb{Z}[x]$ has a $\delta$-structure with the lift of Frobenius given by the formula $F(ax) = ax^{p}$ for any $a\in \mathbb{Z}$.
\end{ex}

\begin{construction}
The forgetful functor $\delta\text{-}\CAlg_{\mathbb{Z},\heartsuit} \rightarrow \CAlg_{\mathbb{Z},\heartsuit}$ has a right adjoint (\cite[Remark 2.7]{BS22}) $\W: \CAlg_{\mathbb{Z},\heartsuit} \rightarrow \delta\text{-}\CAlg_{\mathbb{Z},\heartsuit}$. To give a concrete construction of $\W(A)$ as a $\delta$-ring, notice first that by adjunction, we have an equivalence of sets

$$
\W(A) \simeq \Map_{\CAlg_{\mathbb{Z},\heartsuit}}(\mathbb{Z}[x],\W(A)) \simeq \Map_{\delta-\CAlg_{\mathbb{Z},\heartsuit}}(\mathbb{Z}\{x\},\W(A)) \simeq \Map_{\CAlg_{\mathbb{Z},\heartsuit}}(\mathbb{Z}\{x\}, A) \simeq $$

$$
 \Map_{\CAlg_{\mathbb{Z},\heartsuit}}(\mathbb{Z}[x_{0},x_{1},...], A) \simeq   \Map_{\CAlg_{\mathbb{Z},\heartsuit}}(\bigotimes_{i=0}^{\infty} \mathbb{Z}[x_{i}], A)  \simeq  \prod_{i=0}^{\infty} \Map_{\CAlg_{\mathbb{Z},\heartsuit}}( \mathbb{Z}[x_{i}], A)    \simeq  A^{\mathbb{N}}.
 $$
We can derive the formula for the Frobenius lift endomorphism $F: \W(A) \rightarrow \W(A)$ from this description. Let $a \in \W(A)$ be the image of $x \in \mathbb{Z}\{x\}$ under a specified $\delta$-map $\mathbb{Z}\{x\} \rightarrow \W(A)$, and $(a_{0},a_{1},...)$ be its Witt coordinates where $a_{i}\in A$ are the images of the elements $x_{i}\in \mathbb{Z}\{x\}= \mathbb{Z}[x_{0},x_{1},...] $ under the corresponding map $\mathbb{Z}[x_{0},x_{1},...] \rightarrow A$. Let $(f_{0},f_{1},...):=F(a_{0},a_{1},...)$. From the commutative diagram

$$
\xymatrix{      \mathbb{Z}[x_{0},x_{1},...] \ar[rr]^-{x_{i} \longmapsto x_{i}^{p}+px_{i+1}} \ar[d]&& \mathbb{Z}[x_{0},x_{1},...]  \ar[d]\\
\W(A) \ar[rr]_-{F}&& \W(A)   }
$$
it follows that the coordinates $f_{i}$ are governed by the equations \ref{iterated_equations}. Observe that equations \ref{iterated_equations} imply that iterating $F: \W(A) \rightarrow \W(A)$, we get that $F^{n}(a_{0},a_{1},...) = (f_{0}^{(n)},f_{1}^{(n)},...)$, where the first component $f^{(n)}_{0}$ satisfies

\begin{equation}\label{F_iterated}
f^{(n)}_{0}=a_{0}^{p^{n}}+pa_{1}^{p^{n-1}}+...+p^{n}a_{n}=w_{n}.
\end{equation} 
Following the idea we learnt from the text of J.Borger and B.Wieland \cite{BW05}, we will now obtain the ghost components formulas appearing in Construction \ref{Witt_ghost}. Consider the category $\varphi\text{-}\CAlg_{\mathbb{Z},\heartsuit}$ of commutative rings endowed with an endomorphism $\varphi: A \rightarrow A$. The right adjoint of the forgetful functor $\varphi\text{-}\CAlg_{\mathbb{Z},\heartsuit} \rightarrow \CAlg_{\mathbb{Z},\heartsuit} $ is given by forming the infinite product $A\longmapsto A^{\mathbb{N}}$ endowed with the direct product ring structure, and the endomorphism given by the shift. The Witt vectors $\W(A)$ carries a Frobenius lift endomorphism $F: \W(A) \rightarrow \W(A)$ described above. By adjunction, the projection $\eta: \W(A) \rightarrow A$ gives a map
 
 \begin{equation}\label{ghost_from_plethysm}
 \W(A) \xymatrix{ \ar[rrr]^-{ (\eta, \eta  F, \eta  F^{2},...)  }&&&} A^{\mathbb{N}}
 \end{equation}
Unwinding the definitions and iterating the formula \ref{F_iterated}, we have $\eta  F^{n}(a_{0},a_{1},...)=w_{n}(a_{0},a_{1},...)$, where $w_{n}$ depends only on $a_{0},a_{1},...,a_{n}$ and is given precisely by the formulas \ref{ghost_components}. Using the standard argument, naturality of the ring homomorphism \ref{ghost_from_plethysm} for any $A$, fixes the ring structure on $\W(A)$ uniquely.
 
 \end{construction}
 
 We will give a slightly different approach to the ring structure on Witt vectors which is somewhat more explicit and is based on a concrete linear-algebraic formula for $\W(A)$.

	\subsection{Non-abelian derived functors.}

In this section we will review some basic ideas from the theory of non-abelian derived functors and derived commutative rings. The material of this subsection was initially developed by Brantner-Mathew in \cite{BM19} in the setting over a field, and generalized to more general contexts by Raksit in \cite{R} and by Brantner-Ricardo-Nuiten in \cite{BCN21}. The theory of non-abelian derived functors starts with the notion of the \textbf{sifted completion}. For any small $\infty$-category $\C$, there exists a unique compactly generated $\infty$-category $\mathcal{P}_{\Sigma}(\C)$ which contains all sifted colimits and has the property that for any functor $F: \C \rightarrow \D$ into another $\infty$-category $\D$ containing all sifted colimits, there exists a unique sifted colimit preserving extension $\LF: \mathcal{P}_{\Sigma}(\C) \rightarrow \D$. The functor $\LF$ is often called the \textbf{non-abelian derived functor} of $F$.

\begin{defn}\label{connective_derived_rings}
Let $\CAlg^{\poly}_{\mathbb{Z}}$ be the category of all finitely generated polynomial rings. The $\infty$-category of \textbf{connective derived commutative rings} $\DAlg_{\mathbb{Z},\geq 0}$ is obtained by freely adjoining all sifted colimits to $\CAlg^{\poly}_{\mathbb{Z}}$, $$\DAlg_{\mathbb{Z},\geq 0} :=\mathcal{P}_{\Sigma}(\CAlg^{\poly}_{\mathbb{Z}}).$$ The $\infty$-category of connective derived commutative rings has a forgetful functor $U:\DAlg_{\mathbb{Z},\geq 0} \rightarrow \Mod_{\mathbb{Z},\geq 0}$ to the connective part of the derived $\infty$-category of abelian groups, and the forgetful functor has a left adjoint $\LSym_{\mathbb{Z}}: \Mod_{\mathbb{Z},\geq 0} \rightarrow \DAlg_{\mathbb{Z},\geq 0}$. The functor $\LSym$ is the non-abelian derived functor of the free polynomial ring functor $\Sym_{\mathbb{Z},\heartsuit}: \Mod_{\mathbb{Z}}^{\free,\fg} \rightarrow \CAlg^{\poly}_{\mathbb{Z}}$ defined on the full subcategory $\Mod^{\free,\fg}_{\mathbb{Z}} \subset \Mod_{\mathbb{Z},\geq 0}$ spanned by finitely generated free abelian groups (notice that $\mathcal{P}_{\Sigma}(\Mod_{\mathbb{Z}}^{\free,\fg}) \simeq \Mod_{\mathbb{Z},\geq 0}$). The composition $U \circ\LSym_{\mathbb{Z}}: \Mod_{\mathbb{Z},\geq 0} \rightarrow \Mod_{\mathbb{Z},\geq 0}$ has the structure of a monad, and the $\infty$-category of algebras over this monad is $\DAlg_{\mathbb{Z},\geq 0}$. 
\end{defn}

Definition \ref{connective_derived_rings} generalizes, or "derives" the notion of a commutative ring by allowing it to have non-negative homotopy groups. In fact, it can be shown that the $\infty$-category $\DAlg_{\mathbb{Z},\geq 0} $ is equivalent to the $\infty$-category obtained by inverting weak equivalences in the ordinary category of simplicial commutative rings. But this notion does not capture some objects appearing in nature. For example, derived global sections of the structure sheaf on a non-affine scheme can be endowed with the structure of a cosimplicial, rather than simplicial, algebra. This suggests that it should be possible to derive the notion of a commutative ring in the other direction, i.e. by allowing it to have negative homotopy groups as well. In fact, it turns out that the notion of a commutative ring can be extended in both directions, i.e. by allowing it to have positive as well as negative homotopy groups. To be more precise, it is possible to extend the monad $\LSym_{\mathbb{Z}}: \Mod_{\mathbb{Z},\geq 0} \rightarrow \Mod_{\mathbb{Z},\geq 0}$ to a monad on the whole derived $\infty$-category $\Mod_{\mathbb{Z}}$, so that one gets the free \textbf{derived commutative ring} monad $\LSym_{\mathbb{Z}} : \Mod_{\mathbb{Z}} \rightarrow \Mod_{\mathbb{Z}}$, and the $\infty$-category of algebras for this monad is by definition, the $\infty$-category of \textbf{derived commutative rings}. The paper \cite{R} performs this construction by using the natural filtration on the free symmetric algebra. We will briefly review the Raksit's construction below.

Let $R$ be a discrete ring, and $\Mod_{R}$ the derived $\infty$-category of $R$-module spectra. Then $\Mod_{R}$ has a standard $t$-structure, such that the category $\Mod_{R}^{\free,\fg}$ of finitely generated free $R$-modules freely generates $\Mod_{R,\geq 0}$ by sifted colimits: $\Mod_{R,\geq 0} \simeq \mathcal{P}_{\Sigma}(\Mod_{R}^{\free,\fg})$.

\begin{construction}\label{derived_functors_Cartier}
Let $\T: \Mod_{R}^{\free,\fg} \rightarrow \Mod_{R} $ a functor. We define \textbf{the right-left Kan extension }of $\T$, $\T^{\RL}: \Mod_{R} \rightarrow \Mod_{R}$ as the right Kan extension along the inclusion $\Mod_{R}^{\free,\fg} \subset \Mod_{R,\leq 0}^{\omega}$, followed by the left Kan extension to all of $\Mod_{R}$. 
\end{construction}

\begin{defn}\label{right_left_extendable}

Let $\D$ be an $\infty$-category containing all totalizations and geometric realizations. We say that a functor $G: \Mod_{R,\leq 0}^{\omega} \rightarrow \D$ preserves \textbf{finite coconnective geometric realizations} if for any simplicial object $A_{\ast} $ in $ \Mod_{R,\leq 0}^{\omega}$ which is $n$-skeletal for some $n$ such that the geometric realization $|A_{\ast}|$ belongs to $\Mod_{R,\leq 0}^{\omega}$, the colimit  of $G(A_{\ast})$ exists in $\D$ and the natural map $|G(A_{\ast})| \rightarrow G(|A_{\ast}|)$ is an equivalence. A functor $G: \Mod^{\free,\fg}_{R} \rightarrow \D$ is \textbf{right-left extendable} if the right Kan extension $R G: \Mod_{R,\leq 0}^{\omega} \rightarrow \D$ preserves finite coconnective geometric realizations.

\end{defn}

The following Proposition is the content of \cite[Remark 2.45]{BCN21}.

\begin{prop}\label{rl_equivalence}
Let $\D$ is any $\infty$-category which contains all sifted colimits and finite totalizations. Let $\Fun^{\RL}(\Mod^{\free,\fg}_{R}, \D)$ be the $\infty$-category of right-left extendable functors $F: \Mod^{\free,\fg}_{R} \rightarrow \D$ and $\Fun^{\der}(\Mod_{R},\D)$ be the $\infty$-category of all functors $G: \Mod_{R} \rightarrow \D$ which preserve sifted colimits and finite totalizations of diagrams in $\Mod^{\free,\fg}_{R}$. Then restriction along the inclusion $\Mod^{\free,\fg}_{R} \hookrightarrow \Mod_{R}$ induces an equivalence

$$
\begin{tikzcd}
\Fun^{\der}(\Mod_{R},\D) \ar[r, "\sim"]& \Fun^{\RL}(\Mod_{R}^{\free,\fg},\D).
\end{tikzcd}
$$
\end{prop}

In practice, right-left extendable functors usually arise from additively polynomial functors.

\begin{defn} Let $\C$ be an additive $\infty$-category, and $\D$ an $\infty$-category containing small colimits. An additively polynomial functor $F: \C \rightarrow \D$ of degree $0$ is a constant functor. Assume a polynomial functor of degree $n-1$ has been defined. A functor $F: \C \rightarrow \D$ is \textbf{additively polynomial of degree $n$} if the derivative $DF_{X}(Y) = \fib(F(X\oplus Y) \rightarrow F(X))$ is of degree $n-1$. A functor $F$ is \textbf{additively polynomial} if it is additively polynomial of some degree.

We let $\Fun_{\apoly}(\C_{0} ,\C)$ be the $\infty$-category of additively polynomial functors $\C_{0} \rightarrow \C$.

\end{defn}

\begin{ex}\label{Sym_excisive}
Let $X\in \Mod_{R,\heartsuit}$. The $n$-th symmetric power of $X$ is defined as $$\Sym_{R,\heartsuit}(X) :=(\otimes_{R}^{n} X)_{\Sigma_{n}},$$ where the coinvariants are taken in the naive sense. The functors $\Sym^{\leq n}_{R,\heartsuit} :=\oplus_{i=1}^{n} \Sym^{i}_{R,\heartsuit}: \Mod_{R,\heartsuit} \rightarrow \Mod_{R,\heartsuit}$ are additively polynomial of degree $n$.
\end{ex}

For future reference, we record the following result proven in \cite[Proposition 3.34, Theorem 3.35]{BM19} in the case $R$ is a field, and generalized to an arbitrary derived algebraic context in \cite[Proposition 4.2.14, Proposition 4.2.15]{R}.

\begin{prop}\label{addpoly_and_excpoly}
Let $R$ be a discrete ring, and $\D$ an $\infty$-category containing all limits and colimits. If $F: \Mod^{\free,\fg}_{R} \rightarrow \D$ is an additively polynomial functor of some degree, then it is right-left extendable.
\end{prop}

The right-left extension of an additively polynomial functor of degree $n$ $F: \Mod^{\free,\fg}_{R} \rightarrow \D$ is an \textbf{$n$-excisive} functor (see \cite[Definition 4.2.7]{R} for what it means). The main point here is that $n$-excisive sifted colimit preserving functors from $\Mod_{R}$ to $\D$ correspond to additively polynomial functors $\Mod^{\free,\fg}_{R} \rightarrow \D$ of degree $n$ under the equivalence of Proposition \ref{rl_equivalence}.

\begin{rem}
For any right-left extendable functor $\T: \Mod^{\free,\fg}_{R} \rightarrow \D$, the restriction $\T^{\RL}|_{\Mod_{R,\geq 0}}$ of its right-left extension to the full subcategory of connective objects coincides with the left Kan extension of $F$ along the inclusion $\Mod^{\free,\fg}_{R} \subset \Mod_{R,\geq 0}$. In this case, we will use the notation $\LT: \Mod_{R} \rightarrow \D$ for the derived functor defined on all $\Mod_{R}$.
\end{rem}

\begin{defn}\label{filtered_monad}
Let $\mathbb{Z}_{\geq 0}^{\times}$ be the category of non-negative integers considered with monoidal structure given by multiplication, and $\C$ an $\infty$-category containing all sifted colimits. A \textbf{filtered (sifted colimit preserving) monad} is a lax monoidal functor

 \begin{equation}\label{filtered_monad}
\xymatrix{\T^{\leq \star}: \mathbb{Z}^{\times}_{\geq 0} \ar[r]&      \End_{\Sigma}(\C) . }
\end{equation}

\end{defn}

A.Raksit shows in \cite{R} that for a filtered monad $\T^{\leq \star}$, the colimit $\T:=\underset{\mathbb{Z}_{\geq 0}}{\colim} \;\T^{\leq \star}$ is a monad on $\C$. The advantage of using filtered monads is that many monads are not excisively polynomial on the nose, but have filtrations whose stages are excisively polynomial functors.

\begin{ex}\label{Sym_rl}

The functor

 \begin{equation}\label{filtered_Sym}
\xymatrix{\Sym_{R}^{\leq \star}: \mathbb{Z}^{\times}_{\geq 0} \ar[r]&      \End_{\Sigma}(\Mod_{R})  }
\end{equation}
$$
\xymatrix{ n \ar@{|->}[r] & (X \longmapsto \oplus_{i=1}^{n} \Sym^{i}_{R}(X)) }
$$
is a filtered monad such that $\underset{\mathbb{Z}_{\geq 0}}{\colim} \Sym_{R}^{\leq \star} \simeq \Sym_{R} $. Moreover, the stages of the filtration are sifted colimit preserving excisively polynomial functors. It follows that the functor $\Sym_{R}^{\leq \star}$ lands in the full subcategory $\End_{\Sigma}^{\epoly}(\Mod_{R}) \subset \End(\Mod_{R})), $ which is equivalent to $\Fun^{\apoly}(\Mod^{\free,\fg}_{R}, \Mod_{R}) $ by Proposition \ref{addpoly_and_excpoly}. In other words, the filtered monad $\Sym_{R}^{\leq \star}$ is right-left extended from $\Mod^{\free,\fg}_{R} \subset \Mod_{R}$. As it preserves the connective subcategory $\Mod_{R,\geq 0}\subset \Mod_{R}$, the monad $\Sym_{R}^{\leq \star}$ is in particular, uniquely determined by the filtered monad $\Sym^{\leq \star}_{\Mod_{R,\geq 0}}$ on $\Mod_{R,\geq 0}$.
\end{ex}

\begin{construction}\label{LSym_monad}
Let $\Sym: \Mod^{\free,\fg}_{R} \rightarrow \Mod_{R}$ be the free commutative algebra functor, $\Sym_{R}(M):=\bigoplus_{n\geq 0} (M^{\otimes^{n}})_{\Sigma_{n}}$. Construction \ref{derived_functors_Cartier} supplies a derived functor $\LSym_{R} : \Mod_{R} \rightarrow \Mod_{R}$ which we will refer to as the \textbf{free derived commutative algebra} in $\Mod_{R}$. Then \cite[Construction 4.2.19]{R} shows that the functor $\LSym_{R} \in \End^{\Sigma}(\Mod_{R})$ has a monad structure. This is achieved via a version of Example \ref{Sym_rl} in the derived setting. Let $\End^{\Sigma}_{0}(\Mod_{R,\geq 0})$ be the full subcategory of $\End^{\Sigma}(\Mod_{R})$ consisting of functors which preserve the full subcategory $\Mod^{\free,\fg}_{R} \subset \Mod_{R,\geq 0}$, and $\End^{\Sigma}_{1}(\C_{\geq 0})$ the full subcategory of functors satisfying $\pi_{0}F(X) \in \Mod^{\free,\fg}_{R}$ for any $X\in \C_{0}$. Raksit shows in \cite[Remark 4.2.18]{R} that the inclusion $\End_{0}^{\Sigma}(\Mod_{R,\geq 0}) \subset \End_{1}^{\Sigma}(\Mod_{R,\geq 0})$ admits a \textbf{monoidal} left adjoint $\tau: \End^{\Sigma}_{1}(\Mod_{R,\geq 0}) \rightarrow \End^{\Sigma}_{0}(\Mod_{R,\geq 0})$. As we observed in Example \ref{Sym_rl}, the filtered monad $\Sym^{\leq \star}_{R}$ is uniquely determined by its restriction to the connective subcategory $\Mod_{R,\geq 0} \subset \Mod_{R}$. Moreover, for any $n$, the functor $\Sym^{\leq n}_{R}$ lies in $\End_{1}^{\Sigma}(\Mod_{R})$. We can define a new filtered monad $\LSym^{\leq \star}_{R}$ by the formula $\LSym^{\leq \star}_{\Mod_{R}} := \tau(\Sym_{R}^{\leq \star})$. Moreover, it is an excisively polynomial filtered monad, therefore it extends to a filtered monad on all of $\Mod_{R}$. Defining $$\LSym_{R}:=\underset{\mathbb{Z}_{\geq 0} }{\colim}\;\LSym_{R}^{\leq \star},$$ we get a monad on $\C$. 
\end{construction}

\begin{defn}
A \textbf{derived $R$-algebra} is an algebra over the monad $\LSym_{R}: \Mod_{R} \rightarrow \Mod_{R}$ from Construction \ref{LSym_monad}. We let $\DAlg_{R}:=\Alg_{\LSym_{R}}(\Mod_{R})$ to be the $\infty$-category of derived commutative algebras in $\Mod_{R}$. The forgetful functor $\DAlg_{R} \rightarrow \CAlg_{R}$ commutes with all limits and colimits. Given any other object $A\in \DAlg_{R}$, we define the $\infty$-category of \textbf{derived commutative $R$-algebras} as the under category $$\DAlg_{A}:=(\DAlg_{R})_{A/}.$$
The $\infty$-category $\DAlg_{A}$ is monadic over $\Mod_{A}$ with the monad $\LSym_{A}: \Mod_{A} \rightarrow \Mod_{A}$ which satisfies the formula $$\LSym_{A}(V \otimes_{R} A) \simeq \LSym_{R}(V) \otimes_{R} A$$
for any induced $A$-module $V\otimes_{R} A$.
\end{defn}

In addition to the problem of right-left extending of certain functors on a derived algebraic context, one often encounters a similar problem in the non-linear setting. An example of such a problem is the following. For any $\mathbb{F}_{p}$-algebra $A$, there is a natural Frobenius operation $\varphi: A \rightarrow A$. Can we derive this construction and endow any derived $\mathbb{F}_{p}$-algebra with a Frobenius operation? More generally, assume for concreteness that we have a fixed commutative ring $R$, and a functor $F: \CAlg_{R,\heartsuit} \rightarrow \C$ defined on discrete $R$-algebras and taking values in some presentable $\infty$-category $\C$. Is there a natural way to extend $F$ it to a functor defined on all of $\DAlg_{R}$? This question was tackled by A.Holeman in the text \cite{H22}. We will now briefly review the formalism of non-linear right-left extensions developed by Holeman.

The main idea of Holeman's construction can be summarized as follows. Assume $\C_{0}$ is a discrete additive category, and $\C:=\Fun^{\oplus}(\C_{0},\Mod_{\mathbb{Z}})$. Suppose also that $\T:\C \rightarrow \C$ is a right-left extended monad on $\C$. Let $\Alg_{\T}(\C)^{\poly}$ be the category of polynomial $\T$-algebras, and suppose we are given a functor $F: \Alg_{\T}(\C)^{\poly} \rightarrow \D$ to some presentable $\infty$-category $\D$. Composing $F$ with the free algebra functor $\T: \C_{0} \rightarrow \Alg_{\T}(\C)^{\poly}$, we obtain a functor $F\circ \T: \C_{0} \rightarrow \D$ whose analysis can be done via methods developed in the previous subsection. In particular, there exists a right-left extension $(F\circ \T)^{\RL}: \C \rightarrow \D$. The main idea of Holeman is that under certain conditions, we can endow the functor $F\circ \T$ with an additional structure, keeping which we can upgrade the right-left extension $(F\circ \T)^{\RL}: \C \rightarrow \D$ to a functor $(F\circ \T)^{\RL}: \Alg_{\T}(\C) \rightarrow \D$.
\begin{construction}
Assume $\fcat C$ an $\infty$-category, and $\T: \fcat C\rightarrow \fcat C$ a monad on it. Let

$$\xymatrix{ 
\fcat C\ar@<+1.0ex>[rrr]^-{F} &&& \ar@<+1.0ex>[lll]^-{U}  \Alg_{\T}(\fcat C)
}$$
be the corresponding free-forgetful adjunction. Assume $\fcat D$ is another $\infty$-category. We will be interested in understanding the $\infty$-category of functors $\Fun(\Alg_{\T}(\fcat C), \fcat D)$ in terms of $\Fun(\fcat C, \fcat D)$ endowed with a certain additional data. Since the functors $F$ and $G$ form an adjoint pair, it follows that we have an induced adjunction 

\begin{equation}\label{adjunction_F*}
\xymatrix{ 
\Fun(\Alg_{\T}(\fcat C),\fcat D) \ar@<-1.0ex>[rrr]_-{F^{*}} &&& \ar@<-1.0ex>[lll]_-{U^{*}} \Fun(\fcat C, \fcat D),
}\end{equation}
where the functors $F^*$ and $U^*$ are given by precomposing with $F$ and $U$ respectively, $F^{*}(f):=f\circ F$ and $U^{*}(g):=g\circ U$. Let us denote $\T^{*}$ the resulting monad on $ \Fun(\fcat C, \fcat D)$ so that we have a "forgetful" functor

\begin{equation}\label{T^*-algebras}
\xymatrix{ \Fun(\Alg_{\T}(\fcat C),\fcat D)  \ar[r]& \Alg_{\T^{*}}(\Fun(\fcat C, \fcat D)).      }
\end{equation}
\end{construction}

\begin{rem}
 The $\infty$-category of $\T^{*}$-algebras has an alternative description. The $\infty$-category of functors $\Fun(\fcat C,\fcat D)$ is left-tensored over the monoidal $\infty$-category $\End(\fcat C)$. Given any functor $f: \Alg_{\T}(\fcat C) \rightarrow \fcat D$, the composition $f\circ \F: \fcat C\rightarrow \fcat D$ has the structure of a right module over the monad $\T\in \End(\fcat C)$ acting on $\Fun(\fcat C, \fcat D)$. We have an equivalence of $\infty$-categories
 
 $$
 \Alg_{\T^{*}}(\Fun(\fcat C, \fcat D) \simeq \RMod_{\T}(\Fun(\fcat C, \fcat D)).
 $$
\end{rem}

We will show that the functor $\Fun(\Alg_{\T}(\fcat C),\fcat D)  \rightarrow  \Alg_{\T^{*}}(\Fun(\fcat C, \fcat D)) \simeq \RMod_{\T}(\Fun(\C,\D))$ is close to being an equivalence. More specifically, it becomes an equivalence on the full subcategories of functors which preserve split geometric realizations. We formulate the next proposition for functors preserving sifted colimits, but the proof shows that preservation of split geometric realizations is enough.

\begin{prop}\label{RMod_T}
Assume the monad $\T$ preserves sifted colimits. The adjunction \ref{adjunction_F*} restricts to a monadic adjunction
$$
\xymatrix{ 
\Fun^{\Sigma}(\Alg_{\T}(\fcat C),\fcat D) \ar@<-1.0ex>[rrr]_-{F^{*}} &&& \ar@<-1.0ex>[lll]_-{U^{*}}\Fun^{\Sigma}(\fcat C, \fcat D),
}
$$
and hence we have an equivalence of $\infty$-categories:

$$
\xymatrix{
\Fun(\Alg_{\T}(\fcat C),\fcat D)) \ar[rr]^-{F^{*}}_{\sim}&& \RMod_{\T}(\Fun(\fcat C,\fcat D)) .}
$$
\end{prop}

\begin{proof}
The statement that the adjunction \ref{adjunction_F*} restricts to the full subcategories of functors preserving sifted colimits follows from the assumption that $\T$ preserves sifted colimits. The functor $F^{*}$ clearly commutes with all colimits. For monadicity, it remains to see that it is conservative. Assume $f\rightarrow f$ is a map in $\Fun(\Alg_{\T}(\fcat C), \fcat D)$ which induces an equivalence $f\circ F \simeq g \circ F$. We want to show that for any $A\in \Alg_{\T}(\fcat C)$, the map $f(A) \rightarrow g(A)$ is an equivalence. $A$ has a Bar resolution $A\simeq \colim_{\Delta^{\op}} \T^{\bullet}(A)$ by free $\T$-algebras. Commutation with geometric realizations implies that $$f(A) \simeq \colim_{\Delta^{\op}}f(\T^{\bullet}(A)) \simeq \colim_{\Delta^{\op}} (f\circ F)(U\T^{\bullet -1})(A) \simeq  \colim_{\Delta^{\op}}(g \circ F)(U\T^{\bullet-1}) (A)\simeq \colim_{\Delta^{\op}} g(\T^{\bullet}(A)) \simeq g(A), $$ as desired.
\end{proof}

The next definition gives a necessary condition for the functor $F$ defined on polynomial $\T$-algebras on a derived algebraic context, to be extendable.

\begin{defn}\label{right_left_extendable_nonlinear}
Let $\C_{0}$ be a small additive category, $\C:=\Fun^{\oplus}(\C_{0},\Mod_{\mathbb{Z}})$ and $\T: \C \rightarrow \C$ is a right-left extended monad on $\C$, and $F: \Alg_{\T}(\C)^{\poly} \rightarrow \D$ is a functor to some presentable $\infty$-category $\D$. We say that $F$ is \textbf{right-left extendable} if the functor $F\circ \T: \C_{0} \rightarrow \D$ is right-left extendable in the sense of Definition \ref{right_left_extendable}. Also given a sifted colimit-preseving functor $F: \Alg_{\T}(\C)_{\heartsuit} \rightarrow \D$ from discrete $\T$-algebras in $\C$, we say that is right-left extendable if its restriction to $\Alg_{\T}(\C)^{\poly}$ is.
\end{defn}

The following construction is the content of \cite[Proposition 2.2.15, Construction 2.2.17]{H22}.

\begin{construction}\label{nonlinear_rl_extension}
Assume we are in the situation of Definition \ref{right_left_extendable_nonlinear}. Let $\Fun^{\ext}_{\Sigma}(\C_{\heartsuit},\D)$ be the $\infty$-category of functors $F: \C_{\heartsuit} \rightarrow \D$ which are right-left extendable and preserve sifted colimits in $\C_{\heartsuit}$. The $\infty$-category $\End^{\ext}_{\Sigma}(\C_{\heartsuit})$ has a monoidal structure by \cite[Lemma 2.2.8]{H22}, and the right-left extension functor $\End^{\ext}_{\Sigma}(\C_{\heartsuit}) \rightarrow \End_{\Sigma}(\C)$ is monoidal. The $\infty$-category $\Fun^{\ext}_{\Sigma}(\C_{\heartsuit},\D)$ is right-tensored over $\End^{\ext}_{\Sigma}(\C_{\heartsuit})$, and $\Fun_{\Sigma}(\C,\D)$ is right-tensored over $\End_{\Sigma}(\C)$. The construction of right-left extension $(-)^{\RL}: \Fun^{\ext}_{\Sigma}(\C_{\heartsuit},\D) \rightarrow \Fun_{\Sigma}(\C,\D)$ refines to a lax map of right-tensored $\infty$-categories

$$
\begin{tikzcd}
(-)^{\RL}: \bigl(\Fun^{\ext}_{\Sigma}(\C_{\heartsuit},\D) \bigr)_{\End^{\ext}_{\Sigma}(\C_{\heartsuit})}^{\otimes} \arrow[rr]&& \bigl( \Fun_{\Sigma}(\C,\D)^{\otimes}_{\End_{\Sigma}(\C)}.
\end{tikzcd}
$$
Given that, one can define a functor $\Fun^{\ext}(\Alg_{\T}(\C)^{\poly},\D) \rightarrow \Fun_{\Sigma}(\Alg_{\T}(\C),\D) $ as the composition 

$$
\begin{tikzcd}
\Fun^{\ext}_{\Sigma}(\Alg_{\T_{\heartsuit}}(\C_{\heartsuit}),\D) \arrow[d, swap, "\text{-}\circ \T"] \arrow[rr]&& \Fun_{\Sigma}(\Alg_{\T}(\C),\D) \\
\Mod_{\T_{\heartsuit}}(\Fun^{\ext}_{\Sigma}(\C_{\heartsuit},\D)) \arrow[rr]&& \Mod_{\T}(\Fun_{\Sigma}(\C,\D)) \arrow[u, swap, "\Barr"], 
\end{tikzcd}
$$
\end{construction}

\section{Witt vectors as a $\delta$-Cartier ring.}

 \subsection{Cartier modules.}
 
 In this subsection we will define an $\infty$-category of Cartier modules. Cartier modules package the underlying linear algebraic structure of the ring of Witt vectors.

\begin{defn}
A (discrete) \textbf{Cartier module} is a triple $(M,F,V)$ where $M$ is an abelian group, and $F,V:M\rightarrow M$ a pair of endomorphisms (called \textbf{Frobenius} and \textbf{Verschiebung}) of $M$ satisfying the relation $FV=p$.
\end{defn}

\begin{ex}
 $\mathbb{Z}_{p}$ has a Cartier module structure with $F=\Id$ and $V$ being the operator of multiplication by $p$.
\end{ex}

\begin{ex}
Let $A$ be a commutative ring. Then the ring of $p$-typical Witt vectors $\W(A)$ of $A$ admits a natural Cartier module structure with respect to the Witt vectors Frobenius and Verschiebung.
\end{ex}

\begin{ex}
For a less intuitive example, let $M$ be any $p$-torsion module endowed with an endomorphism $F$. Then $M$ can be considered a Cartier module with $V=0$. 
\end{ex}

\begin{rem}\label{Cartier_Raynaud}
The collection of all Cartier modules form an additive (in fact, abelian) category $\Cart\Mod_{\heartsuit}$. To see this, it is sufficient to observe that the category $\Cart\Mod_{\heartsuit}$ is equivalent to the category of left modules over the associative ring $\mathbb{Z}\langle F,V| FV=p\rangle $ freely generated by two non-commutating elements $F,V$ modulo the relation $FV=p$.

\end{rem}

\begin{construction}\label{(-)[V]}
Let $\varphi\text{-}\Mod_{\mathbb{Z},\heartsuit}$ be the category of abelian groups endowed with an endomorphism. The forgetful functor $\Cart\Mod_{\heartsuit} \rightarrow \varphi\text{-}\Mod_{\mathbb{Z},\heartsuit}$ forgetting the $V$ map and remembering only $F$ has is a left adjoint $(-)[V]: \varphi\text{-}\Mod_{\mathbb{Z},\heartsuit} \rightarrow \Cart\Mod_{\heartsuit}$. For a finitely-generated free abelian group $M$ with an endomorphism $\varphi:M\rightarrow M$, we define $$M[V]=\bigoplus_{n\geq 0} M=\{(x_{1},x_{2},...)|x_{i}\in M\}.$$ The elements of this abelian group are finite sums of the form $\sum_{i=0}^{n} V^{i}(x_{i})$ with Verschiebung defined as a shift, and Frobenius acting by $$F(\sum_{i=0}^{n} V^{i}(x_{i})):=\varphi(x_{0}) + \sum_{i=1}^{n} p V^{i-1}(x_{i}) .$$

\end{construction}

Consider the category $\varphi\text{-}\Mod_{\mathbb{Z},\heartsuit}$ with the symmetric monoidal structure given by the tensor product of underlying $\mathbb{Z}$-modules. In this symmetric monoidal structure, an algebra object is the same as an algebra endowed with an endomorphism. There is a symmetric monoidal structure on $\Cart\Mod_{\heartsuit}$ which makes the functor $(-)[V]$ symmetric monoidal. The tensor product in this symmetric monoidal structure is called the \textbf{box tensor product} of Cartier modules. Below we will follow the exposition of \cite[Chapter 4.2]{AN}.

\begin{defn}
Let $M,N,P$ be Cartier modules. An \textbf{$(F,V)$-bilinear map} $(-,-):M\times N \rightarrow P$ is a bilinear map of abelian groups satisfying the following relations:

$$(F(m),F(n))=F(m,n),$$

$$V(F(m),n)=(m,V(n)), \:\: V(m,F(n))=(V(m),n)$$
for any $m\in M$ and $n\in N$.
\end{defn}

Let us recall the following classical proposition. See \cite[Lemma 4.9]{AN} for more details.

\begin{prop}\label{boxprod}
The functor $\Hom_{(F,V)}(M\times N,-): \Cart\Mod_{\heartsuit} \rightarrow \Set$ is corepresentable by an object $M\boxtimes N$, called the \textbf{box tensor product of Cartier modules}.
\end{prop}

\begin{proof}
We define $M \boxtimes N$ by the formula

$$
M\boxtimes N:=(M\otimes N)[V]/\sim,
$$
where the equivalence relation $\sim$ is additively generated by the relations $(m\otimes Vn)V^{k} \sim (Fm\otimes n)V^{k+1}$ and $(Vm \otimes n)V^{k} \sim (m\otimes Fn)V^{k+1}$ for any $m\in M, n\in N$ and $k\geq 0$.
\end{proof}

\begin{prop}\label{symm_mon_on_CartMod_discrete}
There is a unique symmetric monoidal structure on the category $\Cart\Mod_{\heartsuit}$ which preserves colimits in each variable and makes the functor $(-)[V]: \varphi\text{-}\Mod_{\mathbb{Z},\heartsuit} \rightarrow \Cart\Mod_{\heartsuit}$ symmetric monoidal. The corresponding tensor product on $\Cart\Mod_{\heartsuit}$ is given by the bifunctor

$$ 
\xymatrix{\Cart\Mod_{\heartsuit} \times \Cart\Mod_{\heartsuit} \ar[r]& \Cart\Mod_{\heartsuit}, (M,N)\longmapsto M\boxtimes N, }
$$
where $M\boxtimes N$ is defined by Proposition \ref{boxprod}.
\end{prop}

For the proof of this statement, see \cite[Proposition 4.11]{AN}.

\begin{ex}
The unit object of $\Cart\Mod_{\heartsuit}$ is the submodule $\bigoplus_{n\geq 0}\mathbb{Z} V^{n}([1]) \subset \W(\mathbb{Z})$ of the ring of $p$-typical Witt vectors of $\mathbb{Z}$ consisting of finite Witt series. 
\end{ex}

\begin{defn}
A \textbf{Cartier ring} $(A,F,V)$ is a commutative algebra object in $\Cart\Mod_{\heartsuit}$ with respect to the box tensor product. 
\end{defn}

\begin{ex}\label{A[V]}
Assume $A$ is a commutative ring endowed with an endomorphism $\varphi: A \rightarrow A$. Then $A[V]$ has the structure of a Cartier ring. The multiplication of elements of the form $V^{i}(x)$ and $V^{j}(y)$ for $i,j\geq 1$ and $x, y \in A \subset A[V]$ is governed by the inductive projection formula

\begin{equation}\label{multiplication_of_V}
V^{i}(x)V^{j}(y)=V(V^{i-1}(x))V^{j}(y)=V(V^{i-1}(x)FV^{j}(y))=V(V^{i-1}(x)pV^{j-1}(y))=pV(V^{i-1}(x)V^{j-1}(y)).
\end{equation}
Iterating this formula, and assuming that for instance, $i=j+k$, $k\geq 0$, we obtain:

\begin{equation}\label{formula}
V^{i}(x)V^{j}(y) = p^{j}V^{j}(V^{k}(x)y)=p^{j}V^{j}\bigl( V^{k}(F^{k}(x)y)  \bigr)=p^{j}V^{i}(F^{k}(x)y)
\end{equation}
where we used the projection formula again to compute product $V^{k}(x)y$. Symmetricity of the projection formula implies that the formula \ref{formula} is symmetric in $x,y$. 
\end{ex}

\begin{rem}\label{discrete_V_trivialized}
Note that the multiplication formula \ref{multiplication_of_V} is valid for any Cartier ring $C$. It follows from this formula that for any $C$, and any element of the form $V(x) \in C$, the $p$-th power $V(x)^{p}$ is equal to $0$ modulo $p$. So the $p$-th power Frobenius map $\overline{F}: C/p \rightarrow C/p$ is zero on the image of $\overline{V}: C/p \rightarrow C/p$.
\end{rem}

 \begin{ex}\label{Witt_vectors_Cartier}
 The ring of Witt vectors $\W(A)$ of a discrete commutative ring $A$ has a natural structure of a Cartier ring with respect to Witt vectors Frobenius and Verschiebung.  
 \end{ex}

We now wish to extend the theory of Cartier modules and Cartier rings to the derived setting. The essential idea is the same as before: a derived Cariter module $(M,F,V)$ is an object in the derived category of abelian groups endowed with two endomorphisms $F,V: M \rightarrow M$ and a homotopy $F\circ V \sim p$.

\begin{defn}
Let $\mathbb{N} * \mathbb{N}$ be the free associative monoid on two variables. The functor $\infty$-category $\Mod_{\mathbb{Z}}^{ \B(\mathbb{N} * \mathbb{N})}$ is the $\infty$-category of objects $(M,F,V)$ where $M \in \Mod_{\mathbb{Z}}$, and $F,V: M \rightarrow M$ a pair of endomorphisms. Taking the composition $(M,F,V) \longmapsto (M,F\circ V)$ defines a functor $\Mod_{\mathbb{Z}}^{\B(\mathbb{N}*\mathbb{N})} \rightarrow \Mod_{\mathbb{Z}}^{\B\mathbb{N}} $. There exists another functor $\Mod_{\mathbb{Z}} \rightarrow \Mod_{\mathbb{Z}}^{\B\mathbb{N}}$ endowing any $M$ with the endomorphism $p: M \rightarrow M$. We define the $\infty$-category of Cartier modules $\Cart\Mod$ as the fiber product:

$$
\xymatrix{  \Cart\Mod \ar[d] \ar[r]& \Mod_{\mathbb{Z}}^{\B(\mathbb{N}* \mathbb{N})} \ar[d]\\
\Mod_{\mathbb{Z}} \ar[r]& \Mod_{\mathbb{Z}}^{\B \mathbb{N}}.    }
$$
\end{defn}

\begin{rem}\label{Cartier_ring}
Alternatively, $\Cart\Mod$ is the stable $\infty$-category of left modules over the associative algebra $\mathbb{Z}\langle F,V| FV=p\rangle $:

$$
\Cart\Mod:=\Mod_{\mathbb{Z}\langle F,V| FV=p\rangle }.
$$
\end{rem}

\begin{rem}
There is a $t$-structure on $\Cart\Mod$ whose connective part consists of Cartier modules $M$ whose underlying $\mathbb{Z}$-module is connective. The heart of this $t$-structure is equivalent to $\Cart\Mod_{\heartsuit}$.
\end{rem}

\begin{defn}
Let $M \in \Cart\Mod$. We define $M/V^{n}: =\cofib(V^{n}: M \rightarrow M)$. We say that  $M$ is \textbf{derived $V$-complete} if the natural map $M \rightarrow \underset{\leftarrow}{\lim} \; M/V^{n}$  is a equivalence. We let $\widehat{\Cart}\Mod \subset \Cart\Mod$ be the $\infty$-category of derived $V$-complete Cartier modules. 
\end{defn}

\begin{ex}
Assume $V=0$ in $M$ (or more generally, that is is nilpotent). Then $M/V\simeq M\oplus M[1]$, and the restriction maps $R_{n}: M/V^{n+1} \rightarrow M/V^{n}$ are $\Id \oplus \: 0: M \oplus M[1] \rightarrow  M \oplus M[1]$. The map $M \rightarrow \underset{\leftarrow}{\lim} \:M/V^{n}$ is an equivalence as the limit of the diagram

$$
\xymatrix{ ... \ar[r]^-{0}& M[1] \ar[r]^-{0}& M[1] \ar[r]^-{0}& M[1]   }
$$
is zero.
\end{ex}

Analogously to Construction \ref{(-)[V]} in the discrete case, the forgetful functor $\Cart\Mod \rightarrow \varphi\text{-}\Mod_{\mathbb{Z}} $ has a left adjoint which can be realized as a base change along an algebra map $\mathbb{Z}[F] \rightarrow \mathbb{Z}\langle F,V | FV=p\rangle$. The forgetful functor $\Cart\Mod \rightarrow \varphi\text{-}\Mod_{\mathbb{Z}}$ is the functor of restriction along this map. The left adjoint is given by the base change $$-\underset{\mathbb{Z}[F]}{\otimes} \mathbb{Z}\langle F,V|FV=p\rangle: \varphi\text{-}\Mod_{\mathbb{Z}} \xymatrix{\ar[r]&} \Cart\Mod.$$

Another fundamental construction with Cartier modules is taking the cofiber of $V$ map.

\begin{construction}\label{mod_V}
Given $(M,F,V) \in \Cart\Mod$, let $M/V:=\cofib(V: M \rightarrow M)$. The object $M/V$ carries a lot of additional structure inherited from the Cartier module structure on $M$. The diagram 

$$
\xymatrix{    M \ar[d]_-{=} \ar[r]^-{V}& M \ar[d]^-{F}\\
M \ar[r]_-{p}& M      }
$$
induces a map $\overline{F}: M/V \rightarrow M/p$. Composing it with the further quotient by $V$, we obtain a map $\overline{\varphi}: M/V \rightarrow M/p \rightarrow (M/V)/p. $ 
\end{construction}

\begin{defn}
A $\overline{\varphi}$-module $(X, \overline{\varphi})$ is a pair of $X\in \Mod_{\mathbb{Z}}$ and $\overline{\varphi}: X \rightarrow X/p$. We let $\overline{\varphi}\text{-}\Mod_{\mathbb{Z}}$ be the $\infty$-category of $\overline{\varphi}$-modules.
\end{defn}

\begin{rem}
Construction \ref{mod_V} gives a functor $\Cart\Mod \rightarrow \overline{\varphi}\text{-}\Mod_{\mathbb{Z}}$. It turns out that this functor looses a lot of information and is far from being an equivalence. To get a better hold on the structure of $M/V$ of a Cartier module $M$, one needs to remember morenstructure on $M/V$ than just the $\overline{\varphi}$-map. Consider the functor $(-)/V: \Cart\Mod \rightarrow \Mod_{\mathbb{Z}}$, and let $$E:=\End(-/V: \Cart\Mod \rightarrow \Mod_{\mathbb{Z}})$$ be the algebra of endomorphisms of this functor. Since the functor $(-)/V$ preserves limits and colimits and is conservative on $V$-complete objects, it follows that there is an equivalence

$$
\xymatrix{\widehat{\Cart}\Mod \ar[rr]_-{\sim}^-{(-)/V}&& \Mod_{E},}
$$ 
where the target is the $\infty$-category of left $E$-modules. We will study this algebra in more detail in future work.
\end{rem}

\begin{construction}\label{W_2_linear}
Let $(X, \overline{\varphi}) \in  \overline{\varphi}\text{-}\Mod_{\mathbb{Z}} $. Define $\W_{2}(X)\in \Mod_{\mathbb{Z}}$ as the pull-back

$$
\xymatrix{     \W_{2}(X) \ar[d]_-{F_{1}} \ar[r]^-{R_{1}}& X \ar[d]^-{\overline{\varphi}}\\
X \ar[r]_-{\can}& X/p.     }
$$
This construction depends only on the $\overline{\varphi}$-map on $X$. The data of a lifting of the map $\overline{\varphi}$ to an endomorphism $\varphi: X \rightarrow X$ is equivalent to giving a splitting $s_{1}: X \rightarrow \W_{2}(X)$ of the restriction map $R_{1}: W_{2}(X) \rightarrow X$ as follows:

$$
\xymatrix{  X \ar[rdd]_-{\varphi} \ar@{-->}[rd] \ar[rrd]^-{=} \\
& \W_{2}(X) \ar[d] \ar[r]& X \ar[d]^-{\overline{\varphi}}\\
& X \ar[r]_-{\can} & X/p.     }
$$

\end{construction}

We will now give an alternative concrete construction of the functor $(-)[V]: \varphi\text{-}\Mod_{\mathbb{Z}} \rightarrow \Cart\Mod$ as a filtered colimit of a tower of objects $\W_{n}(X)$, and describe the Cartier module structure on it. The next construction has the advantage that it manifestly sends derived algebra objects to derived algebra objects, which will be important in what follows.

\begin{construction}\label{another_way_of_adding_V}
Let $(X, \varphi) \in \varphi\text{-}\Mod_{\mathbb{Z}}$. We will inductively construct a sequence of objects $\W_{n}(X) \rightarrow ... \rightarrow \W_{2}(X) \rightarrow X$ such that each successive step is defined as the pull-back

\begin{equation}\label{iter_pullback}
\xymatrix{  \W_{i+1}(X) \ar[d]_-{w_{i+1}} \ar[r]^-{R_{i}}& \W_{i}(X) \ar[d]^-{\varphi_{i}}\\
X \ar[r]& X/p^{i}    }
\end{equation}
for a certain canonical map $\varphi_{i}: \W_{i}(X) \rightarrow X/p^{i}$. The maps $\varphi_{i}$ are constructed as follows. First, we let $\varphi_{1}=\overline{\varphi}: X \rightarrow X/p$. Then we define $\varphi_{i+1}: \W_{i+1}(X) \rightarrow X/p^{i+1} $ as the composition

$$
\xymatrix{   \W_{i+1}(X) \ar[r]^-{w_{i}}& X \ar[r]^-{\varphi}& X \ar[r]^-{\can}& X/p^{i+1}    }
$$

The diagram \ref{iter_pullback} implies that the maps $\W_{n+1}(X)\rightarrow \W_{n}(X)$ are extensions by $X$. These extensions are canonically trivialized. Indeed, the diagram

$$
\xymatrix{ \W_{n}(X) \ar@{-->}[rd]^-{s_{n}} \ar[dd]_-{ w_{n}   } \ar@/^2pc/[rrd]^-{\Id}      \\
&   \W_{n+1}(X)  \ar[d]_-{w_{n+1}} \ar[r]^-{R_{n}}& \W_{n}(X) \ar[d]^-{\varphi_{n}}\\
X \ar[r]_-{\varphi} & X  \ar[r]& X/p^{n}    }
$$
provides a splitting of the $n$-th restriction map. In particular, it follows by induction that we have a direct sum splitting $$\W_{n+1}(X) \simeq \underset{i=0}{\bigoplus^{n}} \:XV^{n} ,$$ and an equivalence

$$X[V] \simeq \underset{\rightarrow}{\colim} \: \W_{n}(X),$$
where the colimit is taken over the section maps $s_{n}: \W_{n}(X) \rightarrow \W_{n+1}(X)$ constructed above.  

Now to construct the Cartier module structure on $X[V]$, let $V_{n}: \W_{n}(X) \rightarrow \W_{n+1}(X)$ be the fiber of $R_{n}: \W_{n+1}(X) \rightarrow \W_{n}(X)$. The maps $V_{n}$ commute with restrictions $R_{n}$ and sections $s_{n}$. In addition, there exists a factorization of the map $w_{n+1}: \W_{n}(X) \rightarrow X$ 
$$
\xymatrix{ w_{n+1}: \W_{n+1}(X)\ar[r]^-{F_{n}}& \W_{n}(X) \ar[r]& ... \ar[r]& \W_{2}(X) \ar[r]^-{F_{1}}& X,  }
$$
where the maps $F_{i}: \W_{i+1}(X) \rightarrow \W_{i}(X)$ are defined inductively by the requirement that they commute with restriction maps and each successive map $F_{i}$ is the unique dotted arrow in the diagram below

\begin{equation}\label{Frobenius_finite_length}
\xymatrix{  \W_{i+1}(X) \ar[rr]^-{R_{i}}\ar[rdd]_-{w_{i+1}} \ar@{-->}[rd]^-{F_{i}} && \W_{i}(X) \ar[d]^-{F_{i-1}}  \\
& \W_{i}(X) \ar[d]_-{w_{i}} \ar[r]^-{R_{i-1}}& \W_{i-1}(X) \ar[d]^-{\varphi_{i-1}}\\
& X \ar[r]_-{\can} & X/p^{i-1}.     }
\end{equation}
The maps $F_{n}$ also commute with restrictions and sections. We have an equation up to homotopy

$$
F_{n}V_{n}=p.
$$
The maps $F_{n}$ and $V_{n}$ induce a pair of endomorphisms $F,V: X[V] \rightarrow X[V]$ such that $FV=p$ up to homotopy.
\end{construction}

\begin{prop}
The functor described in Construction \ref{another_way_of_adding_V} is the left adjoint of the forgetful functor.
\end{prop}

\begin{proof}

We will construct unit and counit. The unit map $s: X\rightarrow X[V]$ is just the section of the quotient $X[V] \rightarrow X$. By construction of the Frobenius on $X[V]$, we have a commutative diagram 

$$
\xymatrix{     X \ar[r]^-{\varphi} \ar[d]_-{s} & X \ar[d]^-{s}\\
X[V] \ar[r]_-{F}& X[V],   }
$$
i.e. the map $s$ is a $\varphi$-equivariant map.

Assume now $M$ is a Cartier module. We will construct a Cartier module map $M[V] \rightarrow M$. Notice that for any $n$, there is a pull-back diagram
$$
\xymatrix{   M \ar[d]_-{F^{n}} \ar[r]& M \ar[d]^-{\overline{F^{n}}}\\
M \ar[r]& M/p^{n} . }
$$
By induction, we will construct a sequence of maps $\alpha_{n}: \W_{n}(M) \rightarrow M$ starting with $\alpha_{1}=\Id$ which fit into a commutative diagram 

$$
\xymatrix{   \W_{n}(M) \ar[d]_-{V_{n}} \ar[r]^-{\alpha_{n}}& M \ar[d]^-{V}\\
 \W_{n+1}(M)  \ar[d]_-{F_{n}} \ar[r]^-{\alpha_{n+1}}& M \ar[d]^-{F}\\
 \W_{n}(M) \ar[r]_-{\alpha_{n}}& M    .  }
$$
Given that, there is a commutative diagram 

$$
\xymatrix{     \W_{n}(M) \ar@/_3pc/[ddd]_-{\varphi_{n}}  \ar[d]_-{w_{n}} \ar[r]^-{\alpha_{n}}& M \ar[d]^-{F^{n-1}} \ar[r]& M/V^{n} \ar[ddd]^-{\overline{F^{n}}}  \\
M  \ar[d]_-{F} \ar[r]_-{=}& M \ar[d]^-{F} \\
M \ar[d]_-{\can} \ar[r]_-{=}& M \ar[d]^-{\can} \\
M/p^{n} \ar[r]_-{=}& M/p^{n} \ar[r]_-{=} &M/p^{n} . }
$$
Therefore, we obtain a dotted arrow between two pull-back diagrams:
 $$ \xymatrix{M \ar[dd]_(.7){=} \ar[rr]^(.3){V_{n}...V_{1}}    \ar[rd]_-{=}  && \W_{n+1}(M) \ar[rr] \ar[dd]^(.7){w_{n+1}} \ar@{-->}[dr]^-{\alpha_{n+1}} && \W_{n}(M) \ar[dr] \ar[dd]^(.7){\varphi_{n}}   \\
  &   M \ar[dd]_(.7){=} \ar[rr]^(.3){V^{n}} && M\ar[dd]^(.7){F^{n}} \ar[rr] && M/V^{n} \ar[dd]^(.7){\overline{F^{n}}}  \\
M\ar[rr]_(.7){p^{n}} \ar[rd]_-{=} &&  M \ar[rr]  \ar[rd]^-{=} && M/p^{n} \ar[rd]^-{=}   \\
  &M \ar[rr]_(.7){p^{n}}&& M \ar[rr]&& M/p^{n} , }$$
thus providing the next inductive step of the construction. By taking the colimit of $\alpha_{n}$ maps, we obtain a map of Cartier modules $\alpha: M[V]\simeq \underset{\rightarrow}{\colim}\: \W_{n}(M) \rightarrow M$ which gives the counit of the adjunction.

    \end{proof}

\subsection{Derived Cartier rings.}

The main caveat in setting up the theory of derived Cartier rings as derived algebra objects in $\Cart\Mod$ is that the box product symmetric monoidal structure on $\Cart\Mod_{\heartsuit}$ is not a part of a derived algebraic context, i.e. there is no derived algebraic context structure on $\Cart\Mod$ compatible with the symmetric monoidal structure on the abelian category $ \Cart\Mod_{\heartsuit}$. In fact, the issue already happens for the $\infty$-category $\varphi\text{-}\Mod_{\mathbb{Z}}$. While $\varphi\text{-}\Mod_{\mathbb{Z}}$ has a $t$-structure with a compatible symmetric monoidal structure, this is not a derived algebraic context. The issue is that the unit object $\mathbb{Z} \in \varphi\text{-}\Mod_{\mathbb{Z},\geq 0} $ is not projective. However, in this case it is clear how to define derived algebra objects in $\varphi\text{-}\Mod_{\mathbb{Z}}$: we simply let $\varphi\text{-}\DAlg_{\mathbb{Z}}$ to be the $\infty$-category of derived rings $A$ endowed with an endomorphism $\varphi: A \rightarrow A$. In this definition, $\varphi\text{-}\DAlg_{\mathbb{Z}}$ is monadic over $\varphi\text{-}\Mod_{\mathbb{Z}}$, and the monad is given by forming the $\LSym_{\mathbb{Z}}$-algebra computed on the level of the underlying symmetric monoidal $\infty$-category $\Mod_{\mathbb{Z}}$. We would like to use a similar strategy to define derived algebras internally to the $\infty$-category $\Cart\Mod$. An additional caveat here is that even constructing a symmetric monoidal structure on $\Cart\Mod$ is somewhat subtle, as we can not use concrete formulas for the tensor product as in the case of discrete Cartier modules. Therefore, it is reasonable try to set up the theory of derived Cartier rings without using any symmetric monoidal structure on $\Cart\Mod$ at all, and this can be done as follows.

Recall that a \textbf{quasi-ideal} in a ring $A$ is the data of an $A$-module $I$ endowed with an $A$-linear map $\epsilon: I \rightarrow A$ satisfying the symmetricity property $\epsilon(x)y=x\epsilon(y)$, where the equality is to be hold in $I$ with respect to the $A$-module structure on it. Let $(A,F,V)$ be an algebra object in $\Cart\Mod_{\heartsuit}$, and let $F_{*}A$ be the module obtained by restricting the free $A$-module of rank $1$ along the Frobenius map $F: A \rightarrow A$. For two elements $x\in A$ and $y \in F_{*}A$, we denote $x_{F}*y=y*x_{F}=F(x)y\in F_{*}A$ the result of left (equivalently, right) acting on $y$ by $x$ using the new $A$-module structure. Then the Verschiebung map can be thought as a quasi-ideal $V: F_{*} A \rightarrow A$. Indeed, the $A$-linearity linearity is precisely encoded in the relations $$V(x*y_{F})=V(xF(y))=V(x)y, \:\:V(x_{F}*y)=V(F(x)y)=xV(y)$$ for any $x,y \in A$. Moreover, the symmetricity property is satisfied because

$$
x*V(y)_{F} = xFV(y)=pxy=F(V(x))y = V(x)_{F}*y .
$$
Therefore, we can rephrase the data of an algebra object in the symmetric monoidal category of Cartier modules as the data of an algebra $A$ endowed with an endomorphism $F: A \rightarrow A$ and a quasi-ideal $V: F_{*}A \rightarrow A$ whose composition with the algebra map $F: A \rightarrow F_{*}A$ is the principal ideal $p: F_{*} A \rightarrow F_{*}A$. A similar definition can be given in the derived setting using the theory of derived ideals developed by the author in \cite{pdI}. We will now recall the definition.

\begin{defn}
Let $\Delta^{1}$ be the $1$-simplex considered as a category with two objects $0,1$ and a single non-trivial arrow $0\rightarrow 1$. There is a symmetric monoidal structure $\ast: \Delta^{1}\times \Delta^{1}\rightarrow \Delta^{1} $ on $\Delta^{1}$ defined as:

\begin{equation}\label{join_product}
0\ast 0 = 0, 0\ast 1 = 1, 1\ast 0 = 1, 1\ast 1 = 1
\end{equation}

The functor category $\Fun((\Delta^{1})^{\op},\Mod_{\mathbb{Z},\heartsuit})$ has a \emph{Day convolution symmetric monoidal structure}. The formula for the tensor product of two objects $X^{1}\rightarrow X^{0}$ and $Y^{1}\rightarrow Y^{0}$ is as follows:

$$
(X^{1}\rightarrow X^{0}) \otimes (Y^{1}\rightarrow Y^{0}) \simeq X^{0}\otimes Y^{1} \bigsqcup_{X^{1}\otimes X^{1}} X^{1}\otimes Y^{0} \rightarrow X^{0}\otimes Y^{0}.
$$
We denote  $\Mod_{\mathbb{Z},\heartsuit}^{\Delta^{1}_{\vee}}$ the resulting symmetric monoidal category. A \textbf{quasi-ideal} $(I \rightarrow A)$ is an algebra object in  $\Mod_{\mathbb{Z},\heartsuit}^{\Delta^{1}_{\vee}}$.
\end{defn}

Now consider the stable $\infty$-category $\Mod_{\mathbb{Z}}^{\Delta^{1}_{\vee}}$. It has a Day convolution symmetric monoidal structure, and a $t$-structure, which together form the data of a derived algebraic context. A simple way to construct it is as follows. The main property of the symmetric monoidal $\infty$-category $\Mod_{\mathbb{Z}}^{\Delta^{1}_{\vee}}$ is that the functor $\cofib: \Mod_{\mathbb{Z}}^{\Delta^{1}} \rightarrow \Mod_{\mathbb{Z}}^{\Delta^{1}}$ gives a symmetric monoidal equivalence $\cofib: \Mod_{\mathbb{Z}}^{\Delta^{1}_{\vee}} \simeq \Mod_{\mathbb{Z}}^{\Delta^{1}} $. We can endow the target of this equivalence with the pointwise $t$-structure, and this clearly gives a derived algebraic context structure on $\Mod_{\mathbb{Z}}^{\Delta^{1}}$, which we can transfer via the equivalence to get a derived algebraic context structure on $\Mod_{\mathbb{Z}}^{\Delta^{1}_{\vee}}$.

\begin{defn}
A \textbf{derived ideal} is a derived algebra object in the derived algebraic context  $\Mod_{\mathbb{Z}}^{\Delta^{1}_{\vee}}$. We denote $\DAlg_{\mathbb{Z}}^{\Delta^{1}_{\vee}}$ the $\infty$-category of derived ideals.
\end{defn}

\begin{rem}
Another way to define the $\infty$-category $\DAlg_{\mathbb{Z}}^{\Delta^{1}_{\vee}}$ is by requiring that the functor $\cofib: \Mod_{\mathbb{Z}}^{\Delta^{1}} \rightarrow \Mod_{\mathbb{Z}}^{\Delta^{1}}$ lifts to an equivalence $\cofib: \DAlg_{\mathbb{Z}}^{\Delta^{1}_{\vee}} \simeq \DAlg_{\mathbb{Z}}^{\Delta^{1}} $. In other words, the data of derived ideal $(I \rightarrow A)$ is the same as the data of a derived ring map $A \rightarrow A/I$, where $A/I := \cofib(I\rightarrow A)$.
\end{rem}

\begin{rem}
Another terminology for derived ideals is \textbf{Smith ideals}, since this generalized theory of ideals was initially developed by Jeff Smith in the setting of spectral algebras. 
\end{rem}

\begin{defn}
A \textbf{derived Cartier ring} consists of the following data:
\begin{itemize}
\item A derived ring $A$ together with a derived ring endomorphism $F: A \rightarrow A$,\\
\item A derived ideal $V: F_{*} A \rightarrow A$,\\
\item A homotopy between the composite 

$$
\xymatrix{  F_{*}A \ar[r]^-{V}& A \ar[r]^-{F}& F_{*}A   }
$$
and the principal ideal $p: F_{*}A \rightarrow F_{*}A$.
\end{itemize}
We let $\DCartAlg$ be the $\infty$-category of derived Cartier rings. To define it more precisely starting with the $\infty$-category $\varphi\text{-}\DAlg_{\mathbb{Z}}$, consider the fiber product \begin{equation}\label{pullback_F}
\xymatrix{  \varphi\text{-}\DAlg_{\mathbb{Z}} \underset{\DAlg_{\mathbb{Z}}\Mod }{\times} \DAlg_{\mathbb{Z}}^{\Delta^{1}_{\vee}}  \ar[d] \ar[r]& \DAlg_{\mathbb{Z}}^{\Delta^{1}_{\vee}} \ar[d]\\
\varphi\text{-}\DAlg_{\mathbb{Z}} \ar[r]& \DAlg_{\mathbb{Z}}\Mod ,} \end{equation} where $\DAlg_{\mathbb{Z}}\Mod $ is the $\infty$-category of pairs $(A,M)$ with $A\in \DAlg_{\mathbb{Z}}$ and $M \in \Mod_{A}$, the functor $\varphi\text{-}\DAlg_{\mathbb{Z}} \rightarrow \DAlg_{\mathbb{Z}}\Mod $ sends $(A,F)$ to $(A,F_{*}A)$, and the functor $ \DAlg_{\mathbb{Z}}^{\Delta^{1}_{\vee}} \rightarrow \DAlg_{\mathbb{Z}}\Mod $ sends a derived ideal $(I \rightarrow A)$ to the pair $(A,I)$. In other words, the pull-back \ref{pullback_F} is the $\infty$-category of derived $\varphi$-algebras $(A,F)$ endowed with a derived ideal $V: F_{*}A \rightarrow A$. We now want to impose the condition that the composition $F\circ V$ is homotopic to multiplication with $p$. This is achieved by taking the further pull-back

$$
\xymatrix{  \DCartAlg \ar[d] \ar[r]& \DAlg_{\mathbb{Z}} \ar[d]\\
 \varphi\text{-}\DAlg_{\mathbb{Z}} \underset{\DAlg_{\mathbb{Z}}\Mod }{\times} \DAlg_{\mathbb{Z}}^{\Delta^{1}_{\vee}} \ar[r]& \DAlg_{\mathbb{Z}}^{\Delta^{1}_{\vee}},    }
$$
where the lower horizontal functor sends the data $(A,F,V)$ to the composite derived ideal $F\circ V: A \rightarrow A$ and the right vertical functor endows any derived ring $A$ with the principal derived ideal $p: A \rightarrow A$. We write objects of the $\infty$-category $\DCartAlg$ simply as triples $(A,F,V)$ without directly referencing a prescribed homotopy $F\circ V \sim p$.
\end{defn}

\begin{rem}
It follows from the discussion above that the category $\CartAlg_{\heartsuit}$ of discrete Cartier rings is equivalent to the category $\CAlg(\Cart\Mod_{\heartsuit})$ of commutative algebra objects in the symmetric monoidal category $\Cart\Mod_{\heartsuit}$ with respect to the box tensor product.
\end{rem}

\begin{rem}
Equivalently, the data of a derived Cartier ring can be understood as follows. Suppose $(A,F)$ is a derived ring endowed with an endomorphism. Then lifting $(A,F)$ to a derived Cartier ring is equivalent to giving a derived ring $\overline{A}$, a map $q: A \rightarrow \overline{A}$ and $\overline{F}: \overline{A}\rightarrow F_{*}A/p$ such that there is a pull-back diagram

$$
\xymatrix{  A \ar[r]^-{q} \ar[d]_-{F} & \overline{A} \ar[d]^-{\overline{F}} \\
 F_{*}A \ar[r]_-{\can} & F_{*}A/p  .  }
$$
Indeed, taking horizontal fibers in the diagram above, we obtain the derived ideal $V: F_{*}A \rightarrow A$ and the desired factorization of $p: F_{*}A \rightarrow F_{*}A$.
\end{rem}

\begin{prop}
The forgetful functor $\DCartAlg  \rightarrow \phi\text{-}\DAlg_{\mathbb{Z}}$ is monadic. 

\end{prop}

\begin{proof}
To see conservativity, given a map of derived Cartier rings

$$
\xymatrix{  F_{*}A \ar[r]^-{V} \ar[d]_-{F_{*}f}& A \ar[d]^-{\simeq}_-{f} \ar[r]^-{F} &F_{*}A\ar[d]^-{F_{*}f} \\
F'_{*}A' \ar[r]_-{V'}& A' \ar[r]_-{F'} & F_{*}'A', }
$$
the map $f$ being an equivalence implies that $F_{*}f$ is an equivalence as well, and hence the whole diagram gives an equivalence of derived Cartier rings. Assume $i \longmapsto (A_{i},F_{i},V_{i})$ is a diagram of derived Cartier rings. Let $(A,F)=\lim_{i} (A_{i},F_{i})$. Then taking the limit of $V_{i}$-maps $V_{i}: F_{i,*}A_{i} \rightarrow A_{i}$, we obtain a map $V: F_{*}A \rightarrow A$ giving a derived Cartier ring structure on $A$.
\end{proof}

\begin{rem}
Consequently, the forgetful functor $\DCartAlg \rightarrow \Cart\Mod$ is also monadic.
\end{rem}

\begin{prop}\label{(-)[V]}
The functor $(-)[V] : \varphi\text{-}\Mod_{\mathbb{Z}} \rightarrow \Cart\Mod$ extends to a functor $(-)[V]: \varphi\text{-}\DAlg_{\mathbb{Z}} \rightarrow \Cart\DAlg, $ which is the left adjoint of the forgetful functor.
\end{prop}

\begin{proof}
The left adjoint is given by Construction \ref{another_way_of_adding_V} in the setting of derived rings. If $A\in \varphi \text{-}\DAlg_{\mathbb{Z}}$, we construct a tower $... \rightarrow \W_{n}(A) \rightarrow ...  \rightarrow \W_{2}(A) \rightarrow A$ where each $\W_{n}(A)$ is defined by the iterated pull-back formula \ref{iter_pullback}. By construction, this is a tower of derived rings, the maps $F_{n}: \W_{n+1}(A) \rightarrow \W_{n}(A)$ are derived ring maps and $V_{n}: F_{n,*}\W_{n}(A)  \rightarrow \W_{n+1}(A)$ are derived ideals. Moreover, the sections $s_{n}: \W_{n}(A) \rightarrow \W_{n+1}(A)$ are derived ring maps. Henceforth, $\underset{\rightarrow}{\colim} \: \W_{n}(A) \simeq A[V]$ is a derived Cartier ring. The unit and counit constructed in Construction \ref{another_way_of_adding_V} are derived ring maps, hence this gives the left adjoint of the forgetful functor.
\end{proof}

\subsection{Derived $\delta$-rings.}

In this subsection we will define the $\infty$-category of derived $\delta$-rings. We will give an intrinsic $\infty$-categorical definition of derived $\delta$-rings as derived rings endowed with a derived lift of Frobenius. This definition is clear and concise, however, it is only one of the few natural approaches to setting up this theory. Another way to define derived $\delta$-rings could be via right-left extensions, i.e. starting with the category of ordinary $\delta$-rings, and right-left extending the free $\delta$-ring monad to a monad on $\Mod_{\mathbb{Z}}$. Yet another approach is to use the derived functor of second Witt vectors, and the Witt vectors comonad. We will eventually prove that all these apriori different approaches lead to the same result. The theory of animated $\delta$-rings was developed by Bhatt-Lurie in the foundational text \cite{BL22} on prismatization. The study of derived $\delta$-rings presented below can be considered a non-connective version of Bhatt-Lurie's treatment. 

We begin with the following construction.

\begin{construction}\label{endowing_with_Frobenius}
Let $\DAlg_{\mathbb{F}_{p}[\mathbb{N}]}$ be the $\infty$-category of derived $\mathbb{F}_{p}$-algebras endowed with an endomorphism. The functor $\Frob: \CAlg_{\mathbb{F}_{p}}^{\poly} \rightarrow \CAlg_{\mathbb{F}_{p}[\mathbb{N}], \heartsuit}$ endowing any polynomial $\mathbb{F}_{p}$-algebra $P$ with the $p$-th power Frobenius endomorphism  $\varphi: P \rightarrow P$ left Kan extends to a sifted colimit-preserving functor $\DAlg_{\mathbb{F}_{p},\geq 0} \rightarrow \DAlg_{\mathbb{F}_{p}[\mathbb{N}],\geq 0}.$ The composition $\Frob \circ \Sym_{\mathbb{F}_{p},\heartsuit}: \Mod_{\mathbb{F}_{p}}^{\free,\fg} \rightarrow \DAlg_{\mathbb{F}_{P}[\mathbb{N}]}$ is manifestly filtered by additively polynomial functors. Therefore, by Construction \ref{nonlinear_rl_extension}, this functor is right-left extendable, and there exists a non-linear right-left extension $\Frob:\DAlg_{\mathbb{F}_{p}} \rightarrow \DAlg_{\mathbb{F}_{p}[\mathbb{N}]}$ whose restriction to $\DAlg_{\mathbb{F}_{p},\geq 0}$ coincides with the left Kan extension of the functor endowing any polynomial $\mathbb{F}_{p}$-algebra with the Frobenius map. 
\end{construction}

\begin{rem}
The functor $\DAlg_{\mathbb{F}_{p}}\rightarrow \DAlg_{\mathbb{F}_{p}[\mathbb{N}]} $ endowing a derived commutative $\mathbb{F}_{p}$ with the $p$-th power Frobenius map is monadic and comonadic. The right adjoint is easy to write down for discrete $\mathbb{F}_{p}[\mathbb{N}]$-algebras. Assume $t: A\rightarrow A$ is an algebra with an endomorphism. Then the subalgebra $A^{t=\phi}:=\{a\in A| t(a)=a^{p}\}\subset A$ is universal with respect to all maps of $\mathbb{F}_{p}$-algebras $\eta: B\rightarrow A$ satisfying the property that $\eta(t(b))=(\eta(b))^{p}$ for all $b\in B$. The left adjoint can be described as follows. The universal algebra with respect to maps $\eta:A\rightarrow B$ satisfying $\eta(t(a))=\eta(a)^{p}$ is the quotient of $A$ by the ideal generated by all elements of the form $t(a)-a^{p}$ for $a\in A$. For instance, the object $\Id: \mathbb{F}_{p}[x]\rightarrow \mathbb{F}_{p}[x] $ is sent by the left adjoint to the algebra $\mathbb{F}_{p}[x]/(x-x^{p}) $, and the object $\mathbb{F}_{p}[x]\rightarrow \mathbb{F}_{p}[x], x\longmapsto 0$ goes to $\mathbb{F}_{p}[x]/(x^{p})$.

\end{rem}

Construction \ref{endowing_with_Frobenius} has the following important variant.

\begin{variant}\label{mod_p_variant}
Given any derived ring $A$, the mod $p$ reduction $A/p$ carries a Frobenius endomorphism $A/p \rightarrow A/p$ provided by Construction \ref{endowing_with_Frobenius}. Composing it with the projection $A \rightarrow A/p$, we obtain a map $\overline{\varphi}: A \rightarrow A/p$. This defines a functor $\Frob: \DAlg_{\mathbb{Z}} \rightarrow \overline{\varphi}\text{-}\DAlg_{\mathbb{Z}}$.
\end{variant}

\begin{defn}\label{delta_algebras_1}
Let $\Frob: \DAlg_{\mathbb{F}_{p}} \rightarrow \DAlg_{\mathbb{F}_{p}[\mathbb{N}]}$ be the functor defined in Construction \ref{endowing_with_Frobenius}. We define the $\infty$-category of \textbf{derived $\delta$-rings} as the pull-back

$$
\begin{tikzcd}
\delta\text{-}\DAlg_{\mathbb{Z}} \arrow[r] \arrow[d] & \DAlg_{\mathbb{Z}[\mathbb{N}]}\arrow[d, "-/p"]\\
\DAlg_{\mathbb{F}_{p}} \arrow[r, swap, "\Frob"] & \DAlg_{\mathbb{F}_{p}[\mathbb{N}]}.
\end{tikzcd}
$$
\end{defn}

\begin{rem}\label{variant_definition}
Using the variant of the functor $\Frob$ defined in Remark \ref{mod_p_variant}, equivalently, the $\infty$-category of derived $\delta$-rings can be defined as the pull-back

$$
\begin{tikzcd}
\delta\text{-}\DAlg_{\mathbb{Z}} \arrow[r] \arrow[d] & \DAlg_{\mathbb{Z}[\mathbb{N}]}\arrow[d, "-/p"]\\
\DAlg_{\mathbb{Z}} \arrow[r, swap, "\Frob"] & \overline{\varphi}\text{-}\DAlg_{\mathbb{Z}}.
\end{tikzcd}
$$

\end{rem}

\begin{rem}
Let $\W_{2}(R)$ be defined by the pull-back square as in Construction \ref{W_2_linear} using the derived Frobenius map $\overline{\varphi}: R \rightarrow R/p$. Then the data of a lift of Frobenius is equivalent to the data of a dotted derived ring map

$$
\xymatrix{  R \ar[rdd]_-{\varphi} \ar@{-->}[rd] \ar[rrd]^-{=} \\
& \W_{2}(R) \ar[d] \ar[r]& R \ar[d]^-{\overline{\varphi}}\\
& R \ar[r]_-{\can} & R/p,      }
$$
i.e. a section of the canonical projection $\W_{2}(R) \rightarrow R$. It follows from this that the full subcategory $\delta\text{-}\CAlg_{\mathbb{Z},\heartsuit}$ consisting of discrete rings, is equivalent to the classical category of $\delta$-rings. 
\end{rem}

It follows from the definition of the $\infty$-category $\delta\text{-}\DAlg_{\mathbb{Z}}$ that the forgetful functor $\delta\text{-}\DAlg_{\mathbb{Z}} \rightarrow \DAlg_{\mathbb{Z}}$ commutes with limits and colimits, and therefore admits a left and a right adjoint. We denote the left adjoint \textbf{$\delta$-envelope} functor by $\Env^{\delta}: \DAlg_{\mathbb{Z}} \rightarrow \delta\text{-}\DAlg_{\mathbb{Z}}$, and the right adjoint \textbf{derived Witt vectors ring} functor by $\W: \DAlg_{\mathbb{Z}} \rightarrow \delta\text{-}\DAlg_{\mathbb{Z}}$.

 \begin{construction}\label{Witt_vectors_right_left}
We will give a concrete construction of the right adjoint $\W: \DAlg_{\mathbb{Z}} \rightarrow \delta\text{-}\DAlg_{\mathbb{Z}}$ by using the right-left extension. First, we define the functor $\W : \DAlg_{\mathbb{Z},\leq 0} \rightarrow \delta\text{-}\DAlg_{\mathbb{Z}} $ on coconnective derived rings as the right Kan extension of the functor $\W: \CAlg_{\mathbb{Z},\heartsuit} \rightarrow \delta\text{-}\CAlg_{\mathbb{Z},\heartsuit} $ on discrete rings. Concretely, for any coconnective derived ring $A$ presented as the totalization $A\simeq \Tot(A^{\bullet})$ of a cosimplicial diagram of discrete rings, we have $$\W(A) \simeq \Tor(\W(A^{\bullet})). $$ Endow the ring of Witt vectors of any discrete ring $\W(A)$ with the complete filtration whose stages are $$\W_{\geq n}(A):=\ker(\W(A) \rightarrow \W_{n}(A)).$$ This filtration has the property that associated graded quotients identify with $A$:

$$
\W_{\geq n}(A) /\W_{\geq n+1}(A) \simeq A.
$$
Then by right Kan extension, for any coconnective derived ring $A$, we obtain a complete filtration $\W_{\geq \star}(A)$ on $\W(A)$ such that $$\gr^{\bullet}\W(A) \simeq \bigoplus_{n\geq 0} A.$$
It follows that the functor $\W: \DAlg_{\mathbb{Z}, \leq 0} \rightarrow \DAlg_{\mathbb{Z}, \leq 0}$ preserves coconnective geometric realizations. It follows that $\W: \DAlg_{\mathbb{Z},\leq 0} \rightarrow\delta\text{-}\DAlg_{\mathbb{Z},\leq 0}$ left Kan extends to a functor $\W: \DAlg_{\mathbb{Z}} \rightarrow \delta\text{-}\DAlg_{\mathbb{Z}}$. \end{construction}

We will now give some concrete formulas for the construction of $\W(A)$ of a derived ring. Recall that there is a sequence of functors $\W_{n}: \CAlg_{\mathbb{Z},\heartsuit} \rightarrow \CAlg_{\mathbb{Z},\heartsuit}$, where  for a discrete ring $A$, the ring $\W_{n}(A)$ is defined as the quotient

$$
\W_{n}(A):=\coker(V^{n}: \W(A) \rightarrow \W(A)).
$$
There are multiplicative restriction maps $R_{n}: \W_{n+1}(A) \rightarrow \W_{n}(A)$, multiplicative Frobenii $F_{n}: \W_{n+1}(A) \rightarrow \W_{n}(A)$ and Verschiebung ideals $V_{n}: F_{n,*}\W_{n}(A) \rightarrow \W_{n+1}(A)$. 

\begin{lem}\label{ghosts_lemma}
Let $A$ be a discrete ring. For any $n\geq 1$, the composition 

$$
\xymatrix{  \W_{n+1}(A) \ar[r]^-{F_{n}}& \W_{n}(A) \ar[r]^-{F_{n-1}} & ... \ar[r]^-{F_{2}} & \W_{2}(A) \ar[r]^-{F_{1}}& A   }
$$
is given by the formula
$$
F_{1}\circ ... \circ F_{n}(x_{0},x_{1},...,x_{n}) = x_{0}^{p^{n}} + px_{1}^{p^{n-1}} + ... + p^{n}x_{n}.
$$
\end{lem}

\begin{proof}
We will prove the statement by induction. To proceed with the inductive step, we use the short exact sequence

\begin{equation}\label{fiber_sequence_Verschiebung}
\xymatrix{ \W_{n}(A) \ar[r]^-{V_{n}} & \W_{n+1}(A) \ar[r]& A,    }
\end{equation}
where the map $\W_{n+1}(A) \rightarrow A$ has a set-theoretic Teichmüller representative section $[-]: A \rightarrow \W_{n+1}(A), x \longmapsto (x,0,...,0)$ which commutes with restriction maps and satisfies $F_{n}([x])=[x^{p}]$. Consequently, we have $$F_{1} \circ ... \circ F_{n}([x])=[x^{p^{n}}].$$ Using the fiber sequence \ref{fiber_sequence_Verschiebung}, it remains to observe that 

$$
F_{1}\circ ...  \circ F_{n}\circ V_{n}=F_{1}\circ ... \circ F_{n-1} \circ (F_{n} \circ V_{n}) =pF_{1}\circ ... \circ F_{n-1}.
$$
It follows by writing a general Witt vector $(x_{0},...,x_{n}) \in \W_{n+1}(A)$ as $(x_{0},x_{1},...,x_{n}) = [x_{0}]+ V_{n}(x_{1},...,x_{n})$ that

$$
F_{1}\circ ... \circ F_{n}(x_{0},x_{1},...,x_{n}) = F_{1}\circ ... \circ F_{n}([x_{0}] + V_{n}(x_{1},...,x_{n}))=$$

$$
=F_{1}\circ ... \circ F_{n}([x_{0}]) + F_{1}\circ ... \circ F_{n}\circ V_{n}(x_{1},...,x_{n}) =$$

$$
 =x_{0}^{p^{n}}+pF_{1}\circ ... \circ F_{n-1}(x_{1},...,x_{n})=x_{0}^{p^{n}}+px_{1}^{p^{n-1}}+...+p^{n}x_{n},
$$
as desired.
\end{proof}

\begin{prop}
For any $n\geq 1$, there exists a pull-back diagram in the category of commutative rings

\begin{equation}\label{main_pullback}
\xymatrix{   \W_{n+1}(A) \ar[d]_-{F_{1}\circ ... \circ F_{n}} \ar[r]^-{R_{n}}& \W_{n}(A)  \ar[d]^-{\varphi_{n}}  \\
A \ar[r]_-{\can}& A/p^{n}  , }
\end{equation}
where the right vertical map is given by the formula

$$
\xymatrix{     \varphi_{n}(a_{0},...,a_{n-1})=a_{0}^{p^{n}}+...+p^{n-1}a_{n-1}^{p}   \:\: \mod p^{n}  }
$$
\end{prop}

\begin{proof}
Indeed, it follows by taking horizontal cofibers in the diagram

$$
\xymatrix{     A \ar[d]_-{=}   \ar[rr]^-{V_{n}...V_{1}}&& \W_{n+1}(A)  \ar[d]^-{F_{1} ... F_{n}}\\
A \ar[rr]_-{p^{n}}&& A}
$$
and using Lemma \ref{ghosts_lemma}.
\end{proof}

 \begin{construction}
Let $A$ be a discrete commutative ring. The natural map $\varphi_{n}: \W_{n}(A) \rightarrow A/p^{n}$ for different $n$ satisfy certain compatibilities between each other. The map $\varphi_{1}=\overline{\varphi}: A\rightarrow A/p$ is the $p$-th Frobenius map. Let us denote $\overline{\varphi}^{n}: A \rightarrow A/p$ the $p^n$-th power map sending $a\longmapsto a^{p^n} \mod p$. The diagram 

$$
\xymatrix{ \W(A) \ar[d]_-{F^{n-1}} \ar[r]^-{V}& \W(A) \ar[d]^-{F^{n}} \ar[r]& A \ar[d]^-{\overline{F^{n}}}\\
\W(A) \ar[r]_-{p}& \W(A) \ar[r]& \W(A)/p    }
$$
gives a lifting of the map $\overline{\varphi}^{n}: A \rightarrow A/p$ to a map $\overline{F^{n}}: A \rightarrow \W(A)/p$. The map $\overline{F^{n}}$ fits into a commutative cubical diagram

 $$ \xymatrix {  \W(A) \ar[dd]_(.8){=} \ar[rr]^(.7){=} \ar[rd]^-{V^{n-1}}  && \W(A) \ar[rr] \ar[dd]^(.8){=} \ar[dr]^-{V^{n}} && 0 \ar[dr] \ar[dd]   \\
  &   \W(A) \ar[dd]^(.3){F^{n-1}} \ar[rr]^(.7){V}  && \W(A) \ar[dd]^(.3){F^{n}} \ar[rr]^(.7){R} && A \ar[dd]^(.3){\overline{F^{n}}}  \\
 \W(A) \ar[rr]_(.7){=} \ar[rd]_-{p^{n-1}} &&  \W(A) \ar[rr]  \ar[rd]^-{p^{n}} && 0 \ar[rd]   \\
    &  \W(A) \ar[rr]_(.7){p} && \W(A) \ar[rr]&& \W(A)/p    \\}$$
Taking the cofibers along the diagonal arrows and composing the obtained vertical maps with the further quotient by $V$, we arrive at the commutative diagram of fiber sequences

\begin{equation}\label{map_of_fiber_sequences}
\xymatrix{        \W_{n-1}(A) \ar[rr]^-{V_{n-1}} \ar[d]^-{\varphi_{n-1}} && \W_{n}(A) \ar[d]^-{\varphi_{n}} \ar[rr]&& A \ar[d]^-{\overline{\varphi}^{n}} \\
A/p^{n-1} \ar[rr]_-{p} && A/p^{n}     \ar[rr]&& A/p.  }
\end{equation}
Using right-left extension, we obtain the maps $\varphi_{n}: \W_{n}(A) \rightarrow A/p^{n}$ for any derived ring $A$, which fit into commutative diagrams \ref{map_of_fiber_sequences}.
 
 \end{construction}

 The following proposition is the most crucial step in establishing that the functor $\W$ is the cofree $\delta$-ring functor.

 \begin{prop}\label{delta_splitting}
 Let $A\in \DAlg_{\mathbb{Z}}$. A choice of the data of a lift of Frobenius $\varphi: A \rightarrow A$ is equivalent to the choice of a splitting of $\W(A)$.
 \end{prop}
 
 \begin{proof}
  Assume by induction that for $i\leq n$, the map $\varphi_{i}: \W_{i}(A) \rightarrow A/p^{i}$ factors as a composition $$  \xymatrix{ \varphi_{i}:  \W_{i}(A) \ar[r]^-{w_{i}}& A \ar[r]^-{\varphi}& A \ar[r]^-{\can}& A/p^{i}  } .$$ The case $i=1$ is trivial as $\W_{1}(A)=A$, the map $w_{1}: A\rightarrow A$ is the identity, and then the map $\varphi_{1}=\overline{\varphi}$ satisfies the condition. Given that, $\W_{n+1}(A)$ can be obtained as the iterated pull-back
  
  $$
  \xymatrix{      \W_{n+1}(A) \ar[dddd]_-{w_{n+1}} \ar[r]& \W_{n}(A) \ar[ddd]_-{w_{n}} \ar[r]& ... \ar[dd] \ar[r]& \W_{2}(A) \ar[d]_-{w_{2}} \ar[r]&       A\ar[d]^-{\can \circ \varphi}\\
  &&&A \ar[d]^-{\can \circ \varphi} \ar[r]_-{\can} & A/p   \\
  &&... \ar[d]\ar[r]& A/p^{2}   \\
&A\ar[d]^-{\can \circ \varphi} \ar[r]& ...  \\
A \ar[r]_-{\can}& A/p^{n} }
  $$
  The commutative diagram 
  
    $$
  \xymatrix{      A \ar[rdddd]_-{\varphi^{n}} \ar[rrddd]^-{\varphi^{n-1}}  \ar[rrrdd] \ar[rrrrd]^-{\varphi}\ar[rrrrr]^-{=}&&&&& A\ar[d]^-{\can \circ \varphi}\\
  &&&&A \ar[d]^-{\can \circ \varphi} \ar[r]_-{\can} & A/p   \\
  &&&... \ar[d]\ar[r]& A/p^{2}   \\
&&A\ar[d]^-{\can \circ \varphi} \ar[r]& ...  \\
&A \ar[r]_-{\can}& A/p^{n} }
  $$
then gives a splitting $s: A \rightarrow \W_{n+1}(A)$ of the map $\W_{n+1}(A) \rightarrow A$ in the $\infty$-category of derived rings. We now wish to construct factorization of $\varphi_{n+1}: \W_{n+1}(A) \rightarrow A/p^{n+1}$. We consider the diagram
 $$
 \xymatrix{        \ar@/_2pc/[dd]_(.7){\varphi_{n+1}} \W_{n+1}(A) \ar[d]^-{w_{n+1}} \ar[rr]&& A  \ar@/_2pc/[dd]_(.7){\overline{\varphi}^{n+1}}  \ar[d]^-{\overline{\varphi}^{n+1}}  \ar[rr]&& \W_{n}(A)[1] \ar[d]^-{w_{n}} \ar@/_2pc/[dd]_(.7){\varphi_{n}}  \\
  A \ar[d]^-{\can \circ \varphi} \ar[rr]&& A/p \ar[d]^-{\overline{\varphi}} \ar[rr]&& A[1]   \ar[d]^-{\can \circ \varphi}                          \\
 A/p^{n+1}     \ar[rr]&& A/p \ar[rr] && A/p^{n}[1].}
 $$
Both top and bottom rectangle are commutative diagram exhibiting maps of fiber sequences. By uniqueness of the induced map on fibers, we obtain a homotopy $\varphi_{n+1} \sim \can \circ \varphi \circ w_{n+1}, $ as desired.
 \end{proof}

 \begin{Theor}\label{Witt_right_adjoint}
The functor $\W: \DAlg_{\mathbb{Z}} \rightarrow \delta\text{-}\DAlg_{\mathbb{Z}}$ of Construction \ref{Witt_vectors_right_left} is the right adjoint of the forgetful functor. 
\end{Theor}

\begin{proof}
First, for any derived ring $A$, there is a quotient map $\W(A) \rightarrow A$. Assume that $B$ is a derived $\delta$-ring. Then by Proposition \ref{delta_splitting}, the $\delta$-structure on $B$ gives a unique splitting $B \rightarrow \W(B)$. This gives the desired counit and unit of the adjunction.
\end{proof}

 \begin{construction}\label{LSym_delta}
In the paper \cite{H22} derived $\delta$-rings are initially defined as algebras over a derived monad $\LSym^{\delta}_{\mathbb{Z}}$ acting on $\Mod_{\mathbb{Z}}$, constructed by right-left extending the classical free $\delta$-ring monad. Holeman endows the monad $\Sym^{\delta}_{\mathbb{Z},\heartsuit} $ with the \textbf{$\delta$-filtration} defined as follows (\cite[Definition 2.3.2]{H22}): for the free $\delta$-ring on one variable $\Sym^{\delta}_{\mathbb{Z},\heartsuit} (\mathbb{Z})=\mathbb{Z}\{x\}=\mathbb{Z}[x_{1},x_{2},...]$, we let $\Sym^{\delta,\leq d}_{\mathbb{Z},\heartsuit}$ be the linear span of monomials $x^{j_{n}}_{i_{1}}...x_{i_{n}}^{j_{n}}$ with exponents satisfying the inequality
 
 $$
 \sum_{k=1}^{n} p^{i_{k}}j_{k} \leq d
 $$
 For free $\delta$-rings on several variables, this filtration extends multiplicatively. Then one can show that this filtration yields the structure of a filtered additively polynomial monad on $\Sym^{\delta}_{\mathbb{Z},\heartsuit}$ (\cite[Theorem 2.3.9]{H22}), and therefore is right-left extendable, so that there is a derived monad $\LSym^{\delta}_{\mathbb{Z}}: \Mod_{\mathbb{Z}} \rightarrow \Mod_{\mathbb{Z}}$.

 \end{construction}

\begin{construction}\label{LSym_delta_remark}
 Let us explain why the functor $\LSym^{\delta}: \Mod_{\mathbb{Z}} \rightarrow \Mod_{\mathbb{Z}}$ defined above lifts to a functor $\LSym^{\delta}_{\mathbb{Z}}: \Mod_{\mathbb{Z}} \rightarrow \delta\text{-}\DAlg_{\mathbb{Z}}$ with values in the $\infty$-category of derived rings. First, for any finitely generated free abelian group $M=\mathbb{Z}^{n}$, the polynomial $\delta$-ring $\LSym^{\delta}_{\mathbb{Z}}(M)\simeq \mathbb{Z}\{x_{1},...,x_{n}\}=\bigotimes_{i=1}^{n} \bigotimes_{k\geq 0} \mathbb{Z}[\delta^{k}(x_{i})]$ has an endomorphism defined on generators by the formula
 $$
 \varphi(\delta^{k}x_{i}):=( \delta^{k}x_{i})^{p}+p\delta^{k+1}x_{i}.
 $$
 This endomorphism is clearly a lift of Frobenius, and the construction is natural with respect to maps of polynomial $\delta$-rings. By left Kan extension, we obtain a functor $\LSym^{\delta}_{\mathbb{Z}}: \Mod_{\mathbb{Z},\geq 0} \rightarrow \delta\text{-}\DAlg_{\mathbb{Z}} $. By commutation with sifted colimits, we see that the underlying connective derived ring of $\LSym_{\mathbb{Z}}^{\delta}(M)$ for any $M \in \Mod_{\mathbb{Z},\geq 0}$ is the infinite tensor product
 
 \begin{equation}\label{inf_tensor_product}
 \LSym^{\delta}_{\mathbb{Z}}(M)\simeq \bigotimes_{k\geq 0} \LSym_{\mathbb{Z}}(M) \simeq \LSym_{\mathbb{Z}}(\bigoplus_{k\geq 0}M).
 \end{equation}
 We now observe that any functor satisfying the formula \ref{inf_tensor_product} commutes with finite totalizations. Indeed, it follows from the fact that the functor $M \longmapsto \bigoplus_{k\geq 0}M$ on $\Mod_{\mathbb{Z}}$ commutes with finite totalizations, and the functor $\LSym_{\mathbb{Z}}$ does. Since the forgetful functor $\delta\text{-}\DAlg_{\mathbb{Z}} \rightarrow \DAlg_{\mathbb{Z}}$ commutes with limits, it follows that the functor $\LSym^{\delta}_{\mathbb{Z}}: \Mod_{\mathbb{Z},\geq 0} \rightarrow \delta\text{-}\DAlg_{\mathbb{Z},\geq 0} $ commutes with finite totalizations, and consequently, extends to a sifted colimit and finite totalization preserving functor $\LSym_{\mathbb{Z}}^{\delta}: \Mod_{\mathbb{Z}} \rightarrow \delta\text{-}\DAlg_{\mathbb{Z}}$ lifting the functor of Construction \ref{LSym_delta}. 
\end{construction}

 In the rest of this subsection, we will prove that Definition \ref{delta_algebras_1} of derived $\delta$-rings is equivalent to algebras over the monad $\LSym^{\delta}_{\mathbb{Z}}$ defined in Construction \ref{LSym_delta}. This is the subject of \cite[Theorem 2.4.4]{H22}, and the sketch of the proof is given in Holeman's paper. We will give a detailed argument here. 
 
 \begin{prop}
 There is an equivalence of $\infty$-categories
 
 $$
 \delta\text{-}\DAlg_{\mathbb{Z}} \simeq \Alg_{\LSym^{\delta}}(\Mod_{\mathbb{Z}}) \simeq \CoAlg_{\W}(\DAlg_{\mathbb{Z}}).
 $$
 \end{prop}
 
 \begin{proof}
The equivalence $ \delta\text{-}\DAlg_{\mathbb{Z}}\simeq \CoAlg_{\W}(\DAlg_{\mathbb{Z}})  $ follows from comonadicity of the forgetful functor and Theorem \ref{Witt_right_adjoint}. To prove the equivalence $ \delta\text{-}\DAlg_{\mathbb{Z}} \simeq \Alg_{\LSym^{\delta}}(\Mod_{\mathbb{Z}})$, it suffices to show that the left adjoint of the forgetful functor is the $\LSym^{\delta}$-functor given by right-left extension as in Construction \ref{LSym_delta}. Construction \ref{LSym_delta_remark} provides a natural map $\Env^{\delta}\circ \: \LSym(M) \rightarrow \LSym^{\delta}(M) $ for any $M$ and using commutation with filtered colimits, it is enough to prove it is an equivalence in the case $M \in \Mod^{\omega}_{\mathbb{Z}}$. Then we want to show that there exists a natural equivalence of mapping spaces

\begin{equation}\label{map_maps}
\Map_{\DAlg_{\mathbb{Z}}}(\LSym^{\delta}(M), B) \simeq \Map_{\DAlg_{\mathbb{Z}}}(\Env^{\delta}(\LSym(M)), B)
\end{equation}
for any derived algebra $B$. The right hand-side can be identified with the mapping space $\Map_{\mathbb{Z}}(M,\W(B))$ in the $\infty$-category $\Mod_{\mathbb{Z}}$, and using the $V$-adic filtration of $\W(B)$, this mapping space can be identified with the infinite product $\prod_{n\geq 0}\Map_{\mathbb{Z}}(M,B)$, so the functor $B \longmapsto \Map_{\mathbb{Z}}(M,B) $ commutes with geometric realizations and totalizations of derived rings. To compute the left-hand side, we use the fact that $\LSym^{\delta}(M) \simeq \otimes_{i \geq 0} \LSym(M)$ as a derived algebra, and hence the corepresentable functor $\Map_{\DAlg_{\mathbb{Z}}}(\LSym^{\delta}(M), -)$ also identifies with the infinite product of mapping spaces $\Map_{\mathbb{Z}}(M,B)$, and hence commutes with geometric realizations and totalizations in the variable $B$. Presenting any derived ring as a simplicial-cosimplicial ring as in \cite[Corollary 5.29]{BCN21}, it is now enough to show the equivalence \ref{map_maps} for a discrete commutative ring $B$, in which case the statement reduces to the classical one.

 \end{proof}

We finish this subsection with some concrete computations of the derived Frobenius and the derived Witt vectors in the special case of $p$-torsion rings. For an $\mathbb{F}_{p}$-algebra $A$, the derived Frobenius takes the form $A \rightarrow A/p \simeq A\oplus A[1]$, where the target $A/p \simeq A\oplus A[1]$ is a trivial square-zero extension. It can be expressed in terms of the Frobenius endomorphism and the Bockstein differential, as Proposition \ref{composition} below shows.  Let us recall the following construction first.

\begin{construction}\label{using_Tate_valued}
Recall from \cite{NS} that any $\mathbb{E}_{\infty}$-ring spectrum $A$ carries a \textbf{Tate-valued Frobenius} map $A\rightarrow A^{t C_{p}}$. where the target is the Tate cohomology spectrum of $A$ with respect to the trivial $C_{p}$-action. If $A$ is connective, then the map $A\rightarrow A^{t C_{p}}$ factors through the connective cover of the target $A\rightarrow \tau_{\geq 0}A^{t C_{p}}$. Assume $A$ is a discrete ring. Then $\tau_{\geq 0}A^{t C_{p}}$ is an $\mathbb{E}_{\infty}$-algebra over $\tau_{\geq 0}\mathbb{Z}^{t C_{p}}$, whose homotopy is the polynomial ring $\mathbb{F}_{p}[\sigma]$, and we have an identification $\tau_{\geq 0}A^{t C_{p}}/\sigma = \tau_{\geq 0}A^{t C_{p}}\otimes_{\tau_{\geq 0}\mathbb{Z}_{p}^{t C_{p}}} \mathbb{F}_{p} \simeq A/p$. Composing the map $A \rightarrow \tau_{\geq 0}A^{t C_{p}}$ with the quotient $\tau_{\geq 0}A^{t C_{p}} \rightarrow \tau_{\geq 0}A^{t C_{p}} /\sigma $, we obtain a map $A \rightarrow A/p$, which is equivalent to the derived Frobenius $\overline{\varphi}: A\rightarrow A/p$ constructed earlier using left Kan extension. Indeed, to check that they are equivalent for a general discrete ring $A$, it suffices to check it for a $p$-torsion free $A$, in which case $A/p$ is discrete, and the map $ A \rightarrow \tau_{\geq 0}A^{t C_{p}} \rightarrow \tau_{\geq 0}A^{t C_{p}}/\sigma \simeq A/p$ is the same as the composition $A\rightarrow \tau_{\geq 0}A^{t C_{p}} \rightarrow \pi_{0} \tau_{\geq 0} A^{t C_{p}} $ which is easily seen to be the $p$-th power map.
\end{construction}

\begin{ex}\label{Z/p^2}
Consider the case $R=\mathbb{F}_{p}$. Let us identify the derived Frobenius operation $\overline{\varphi}: \mathbb{F}_{p} \rightarrow \mathbb{F}_{p} \otimes_{\mathbb{Z}} \mathbb{F}_{p}\simeq \mathbb{F}_{p}\oplus \mathbb{F}_{p}[1]$. Let $\beta_{2}: \mathbb{F}_{p} \rightarrow \mathbb{F}_{p}[1]$ be the Bockstein differential. Since $\mathbb{F}_{p}\oplus \mathbb{F}_{p}[1] \simeq \tau_{\geq 0}\mathbb{F}_{p}^{t C_{p}}/\sigma$, it follows from \cite[Proposition IV.1.15]{NS} that the derived Frobenius on $\mathbb{F}_{p}$ is exprerssed via the Bockstein differential as $(\Id,\beta_{2}):\mathbb{F}_{p} \rightarrow \mathbb{F}_{p} \oplus \mathbb{F}_{p}[1]$. The pull-back diagram

\begin{equation}\label{Z_p^2}
\xymatrix{ \mathbb{Z}/p^{2} \ar[d] \ar[r]& \mathbb{F}_{p} \ar[d]^-{(\Id,\beta_{2})}\\
\mathbb{F}_{p} \ar[r]_-{(\Id,0)}& \mathbb{F} \oplus \mathbb{F}_{p}[1],    }
\end{equation}
recovers the equivalence $\W_{2}(\mathbb{F}_{p})\simeq \mathbb{Z}/p^{2}$ as the square-zero extension of $\mathbb{F}_{p}$ via the Bockstein differential. 
\end{ex}

More generally, assume $A$ is a $p$-torsion free ring, and we want to compute the $\overline{\varphi}$-operation on $A/p$. We have already seen in Example \ref{Z/p^2} that the $\overline{\varphi}$-operation on $\mathbb{F}_{p}$ is given by $(\Id, \beta_{2}): \mathbb{F}_{p} \rightarrow \mathbb{F}_{p}\oplus \mathbb{F}_{p}[1]$. Since $A/p$ comes as the base change of the derived ring $A$, it also carries the Bockstein differential map $\beta_{2}: A/p \rightarrow A/p[1]$. Let us denote the Frobenius endomorphism by $\varphi: A/p \rightarrow A/p$. 

\begin{prop}\label{composition}
The map $\overline{\varphi}: A/p \rightarrow (A/p)/p \simeq A/p \oplus A/p[1]$ is the composition 

$$
\xymatrix{   A/p \ar[r]^-{\varphi}& A/p \ar[r]^-{(\Id, \beta_{2})} & A/p \oplus A/p[1].     }
$$
\end{prop}

\begin{proof}
Using Construction \ref{using_Tate_valued}, the map $A/p \rightarrow A/p \oplus A/p[1]$ is the composition $$\xymatrix{   A/p \ar[r]& \tau_{\geq 0} (A/p)^{t C_{p}} \ar[r]^-{\sigma \longmapsto 0} & A/p \oplus A/p[1]     } .$$ Using the proof of \cite[Proposition IV.1.16]{NS}, the Tate valued Frobenius factors as a map of spaces into the composition
$$
\xymatrix{  A/p \ar[r]^-{\Delta}& \bigl((A/p)^{\times^{p}}\bigr)^{\Sigma_{p}} \ar[r]^-{m^{\Sigma_{p}}}&   (A/p)^{\Sigma_{p}}   \ar[r]^-{\simeq } & (A/p)^{C_{p}} \ar[r]& \tau_{\geq 0}(A/p)^{t C_{p}}  ,}
$$
where $\Delta$ is the set-theoretic diagonal, the map $m^{\Sigma_{p}}$ is induced on $\Sigma_{p}$-fixed points by the multiplication, and the last two map is the natural map from fixed points to the Tate cohomology. Note that $ (A/p)^{C_{p}} \simeq A/p $, and the composition of the first two maps is the Frobenius endomorphism $\varphi: A/p \rightarrow A/p$. It remains to identify the composition $\xymatrix{A/p \ar[r]& \tau_{\geq 0}(A/p)^{tC_{p}}\ar[r]^-{\sigma\longmapsto 0}& A/p \oplus A/p[1]  }$ with the map $(\Id,\beta_{2})$. It follows from the initial case $A/p=\mathbb{F}_{p}$. 
\end{proof}

\begin{warning}
Notice that the order of operations is important. The composition \begin{equation}\label{another_comp}\xymatrix{ A/p \ar[r]^-{(\Id, \beta_{2})} &  A/p \oplus A/p[1] \ar[r]^-{(\varphi, \varphi[1])}& A/p \oplus A/p[1]  } \end{equation} is the natural map obtained by taking the quotient by $p$ of the map $\overline{\varphi}_{A}: A \rightarrow A/p$. This map is \emph{not} equivalent to $\overline{\varphi}_{A/p}$ in general. In fact, a homotopy between the compositions $\beta_{2}\varphi$ and $\varphi[1]\beta_{2}$ is equivalent to a lift of Frobenius modulo $p^{2}$. \end{warning}

\begin{prop}
Let $A$ be a derived $\mathbb{F}_{p}$-algebra. Then each map $\W_{n+1}(A) \rightarrow \W_{n}(A)$ is a square-zero extension.
\end{prop}

\begin{proof}
Indeed, since $p=0$ in $A$, it follows that $A/p^{n} \simeq A\oplus A[1]$ and $\W_{n+1}(A)$ is a pull-back

$$
\xymatrix{      \W_{n+1}(A) \ar[d]_-{w_{n+1}} \ar[r]^-{R_{n}}& \W_{n}(A) \ar[d]^-{\varphi_{n}}\\
w_{n+1,*}A \ar[r]_-{(\Id,0)}& w_{n+1,*}(A\oplus A[1]).    }
$$
The fiber of the map $R_{n}$ then identifies with the fiber of the lower horizontal map, so that $$ \fib(R_{n}) \simeq w_{n+1,*}A$$ as a trivial non-unital $\W_{n+1}(A)$ algebra, as desired.
\end{proof}

\subsection{Derived $\delta$-Cartier rings.}

\begin{defn}

A \textbf{pre-$\delta$-Cartier ring} $(A,F,V,\delta)$ is the following data:

\begin{itemize}

\item A commutative ring $A$;

\item A $\delta$-structure with a lift of Frobenius $F: A \rightarrow A$;

\item A quasi-ideal $V: F_* A \rightarrow A$ such that $$FV=p$$.

\end{itemize}

\end{defn}

\begin{rem}
Consider the equation $FV(a)=pa$ for any $a\in A$. Unwinding the definition, we obtain the following equation:

$$
V(a)^{p} +p\delta(Va)=pa.
$$
Using projection formula, we have $V(a)^p = p^{p-1}V(a^p)$. If the ring $A$ is $p$-torsion free, we can solve this uniquely for $\delta(Va)$.

\begin{equation}\label{delta_V}
\delta(Va)=a-p^{p-2}V(a^{p})
\end{equation}
 However, if there is $p$-torsion, there is an ambiguity in terms of what $\delta(Va)$ can be, as the following example shows.

\end{rem}

\begin{ex}
Consider a square-zero extension ring of the form $\mathbb{Z}\oplus I$, where $I$ is an $\mathbb{F}_{p}$-vector space. Let $V:= p \oplus \Id_I$ and $F=\Id_{\mathbb{Z}} \oplus \: 0_{I}$. Then $F$ and $V$ endow $A$ with the structure of a Cartier ring. Letting $\delta(n,x)=(\delta(n), 0)$, where $\delta(n)=\frac{n-n^{p}}{p}$, we promote it to the structure of a pre-$\delta$-Cartier ring. However, $\delta(Vx)=\delta(x)=0\neq x$ if $x\in I$ is non-zero, hence the formula above does not hold.
\end{ex}

The notion of a $\delta$-Cartier ring is obtained by enforcing the relation \ref{delta_V} in the notion of a pre-$\delta$-Cartier ring.

\begin{defn}\label{delta_Cartier}

A \textbf{ $\delta$-Cartier ring} $C$ is a pre-$\delta$-Cartier ring $A$ for which the formula \ref{delta_V} holds for any $a\in A$.

\end{defn}

Let $\delta\text{-}\CartCAlg$ be the category of $\delta$-Cartier rings. There is a forgetful functor $\delta\text{-}\CartCAlg \rightarrow \delta\text{-}\CAlg_{\mathbb{Z}}$. We will show below that the the left adjoint of this functor is given by freely adjoining $V$ to a $\delta$-ring.

\begin{construction}\label{(-)[V]_for_delta}
Assume $A$ is a $p$-torsion free $\delta$-ring. Let us describe a concrete $\delta$-Cartier ring structure on $A[V]$. We want to construct a $\delta$-operation $\delta:A[V] \rightarrow A[V]$ such that $F(a)=a^{p}+p\delta(a)$. Elements of $A[V]$ are finite sums $\sum_{n\geq 0}V^{n}(a_{n})$, and for any $a\in A[V]$, we have $FV^{n}(a)=pV^{n-1}(a)$ in $A[V]$. Define 

\begin{equation}\label{delta_formula}
\delta(V^{n}(a))=V^{n-1}(a)-p^{(p-1)n-1}V^{n}(a^{p}),
\end{equation}
for any $a\in A$. Extend $\delta$ to arbitrary finite sums $\sum_{n\geq 0}V^{n}(a_{n})$ by the $\delta$-ring sum formula.

\end{construction}

\begin{prop}

 The map $\delta: A[V] \rightarrow A[V]$ gives a $\delta$-Cartier ring structure on $A[V]$.

\end{prop}

\begin{proof}

Notice that the formula \ref{delta_V} for elements of the form $a\in A \subset A[V]$ is a special case of  \ref{delta_formula} for $n=1$. Moreover, since we have

$$
p^{(p-1)n-1}V^{n}(a^{p}) = p^{p-2}p^{(p-1)(n-1)}V(V^{n-1}(a^{p}))$$

$$
=p^{p-2} V\bigl( p^{(p-1)(n-1)}V^{n-1}(a^{p})\bigr)
$$

$$
=p^{p-2}V \bigl( V^{n-1}(a)\bigr)^{p},
$$
it follows that the relation \ref{delta_V} in fact holds for any $a\in A[V]$. Moreover, we we have
 \begin{equation}\label{equation_V_delta}
 V^{n}(a)^{p}+p\delta(V^{n}(a)) = p^{(p-1)n}V^{n}(a^{p}) + p\bigl(    V^{n-1}(a)-p^{(p-1)n-1}V^{n}(a^{p})  \bigr) = pV^{n-1}(a),
 \end{equation} 
hence the Cartier ring condition $F(V^{n}(a))=(V^{n}(a))^{p}+p\delta(V^{n}(a))=pV^{n-1}(a)$ holds as desired. Hence we have $F(a)=a^{p}+p\delta(a)$  for any $a\in A[V]$. Since $F$ is a ring homomorphism and $A[V]$ is $p$-torsion free, it follows that $\delta: A[V] \rightarrow A[V]$ satisfies the $\delta$-ring axioms and the relation \ref{delta_formula} holds by definition. Therefore, $A[V]$ is a $\delta$-Cartier ring.
\end{proof}

\begin{prop}
Assume $A$ is a $p$-torsion free $\delta$-ring and $B$ is a $\delta$-Cartier ring. There is a bijection

$$
\Hom_{\delta\text{-}\CAlg_{\mathbb{Z}}}(A,B) \simeq \Hom_{\delta\text{-}\CartCAlg}(A[V],B),
$$
where $A[V]$ is as described in Construction \ref{(-)[V]_for_delta}.
\end{prop}

\begin{proof}
Indeed, since $A[V]$ is the free Cartier ring on $A$, any $\delta$-map $A \rightarrow B$ uniquely extends to a Cartier ring map $A[V] \rightarrow B$. Since the relation \ref{delta_V} holds in $A[V]$ and $B$, it follows that the map $A[V] \rightarrow B$ necessarily commutes with $\delta$-structures, hence is a map of $\delta$-Cartier rings. 
\end{proof}

\begin{defn}
Let $\delta\text{-}\CartCAlg^{\poly} \subset \delta\text{-}\CartCAlg$ be the subcategory of $\delta$-Cartier rings of the form $\mathbb{Z}\{x_{1},...,x_{n}\}[V]$ for all $n\geq 1$, and 

$$
\delta\text{-}\DCartAlg_{\geq 0}:=\mathcal{P}_{\Sigma}( \delta\text{-}\CartCAlg^{\poly}).
$$
Then $\delta\text{-}\DCartAlg_{\geq 0}$ is monadic over $\Mod_{\mathbb{Z},\geq 0}$ with the monad given by $$\LSym^{\delta}[V](M):=\LSym^{\delta}(M)[V]:=\colim_{i\in I} \mathbb{Z}\{ x_{1},...,x_{\rk(F_{i})}\}[V]$$
for any presentation of $M$ as a sifted colimit $M=\colim_{i\in I} F_{i}$ of $F_{i}\in \Mod_{\mathbb{Z}}^{\free,\fg}$. This monad extends to a sifted colimit preserving monad $\LSym^{\delta}[V]: \Mod_{\mathbb{Z}} \rightarrow \Mod_{\mathbb{Z}}$. Now we let the $\infty$-category of \textbf{derived $\delta$-Cartier rings} to be

$$
\delta\text{-}\DCartAlg:=\Alg_{\LSym^{\delta}[V]}(\Mod_{\mathbb{Z}}).
$$
Any derived $\delta$-Cartier ring has an underlying Cartier module. We say that a derived $\delta$-Cartier ring $A$ is \textbf{$V$-complete} if the underlying Cartier module is $V$-complete.
\end{defn}

The forgetful functor $\delta\text{-}\DCartAlg \rightarrow  \delta\text{-}\DAlg_{\mathbb{Z}}$ is monadic. By construction, the left adjoint $(-)[V]: \delta\text{-}\DAlg_{\mathbb{Z}} \rightarrow \delta\text{-}\DCartAlg$ fits into the commutative diagram

\begin{equation}\label{add_V}
\xymatrix{
\delta\text{-}\DAlg_{\mathbb{Z}} \ar[d] \ar[r]^-{(-)[V]} & \delta\text{-}\DCartAlg \ar[d]\\
\varphi\text{-}\DAlg_{\mathbb{Z}} \ar[r]_-{(-)[V]}&  \DCartAlg.
}
\end{equation}
commutes.

It is now formal to construct a functor $\W: \DAlg_{\mathbb{Z}} \rightarrow \delta\text{-}\widehat{\Cart}\DAlg $ such that the composition $\DAlg_{\mathbb{Z}} \rightarrow \delta\text{-}\widehat{\Cart}\DAlg \rightarrow \delta\text{-}\DAlg_{\mathbb{Z}} $ is the derived Witt vectors functor.

\begin{construction}\label{two_different_Witt_vectors}

The functor $(-)[[V]]:\delta\text{-}\DAlg_{\mathbb{Z}} \rightarrow \delta\text{-}\widehat{\Cart}\DAlg$  identifies with the left adjoint of the forgetful functor, and the diagram

$$
\xymatrix{\delta\text{-}\DAlg_{\mathbb{Z}} \ar[rd] \ar[rr]^-{(-)[[V]]}&& \delta\text{-}\widehat{\Cart}\DAlg    \ar[ld]^-{(-)/V} \\
& \DAlg_{\mathbb{Z}} .}
$$
commutes. Passing to right adjoints, we obtain a functor $\W: \DAlg_{\mathbb{Z}} \rightarrow \delta\text{-}\widehat{\Cart}\DAlg$ which fits into a commutative triangle

$$
\xymatrix{   \delta\text{-}\DAlg_{\mathbb{Z}} &&\ar[ll]_-{\forget_{V}}  \delta\text{-}\widehat{\Cart}\DAlg\\
& \ar[lu]^-{\W} \DAlg_{\mathbb{Z}} \ar[ru]_-{\W},  }
$$
where the upper horizontal functor is the forgetful one. 
\end{construction}

The main reason for introducing the $\infty$-category $\delta\text{-}\widehat{\Cart}\DAlg$, is due to the following Theorem.

\begin{Theor}\label{Witt_equivalence_1}
The adjunction of Construction \ref{two_different_Witt_vectors} gives an equivalence

$$
\xymatrix{(-)/V : \delta\text{-}\widehat{\Cart}\DAlg   \ar@<+1.0ex>[rrr]_-{\sim} &&& \ar@<+1.0ex>[lll] \DAlg_{\mathbb{Z}} : \W, }
$$
where the lift of the functor $\W$ is as defined in Construction \ref{two_different_Witt_vectors}.
\end{Theor}

\begin{proof}
Indeed, for any derived ring $R$, the counit $\W(R)/V \simeq R$ is an equivalence. Similarly, for any derived $\delta$-Cartier ring $A$, the unit of the adjunction $f: A\rightarrow \W(A/V)$ is induced by the quotient map $A\rightarrow A/V$. The map $f$ is therefore an equivalence on $V$-reductions, hence an equivalence of derived $\delta$-Cartier rings by $V$-completeness.
\end{proof}

As an immedate consequence of symmetric monoidality, we also have the next relative version. 

\begin{cor}\label{Witt_equivalence_relative}
Let $A\in \DAlg_{\mathbb{Z}}$ endowed with the derived Frobenius operation $A \rightarrow A/p$. Then there is an equivalence of $\infty$-categories

$$
\xymatrix{(-)/V : \delta\text{-}\widehat{\Cart}\DAlg_{\W(A)}  \ar@<+1.0ex>[rrr]_-{\sim} &&& \ar@<+1.0ex>[lll] \DAlg_{A} : \W, }
$$
\end{cor}

\begin{rem}
Assume that $C$ is a $p,V$-torsion free $V$-complete $\delta$-Cartier ring whose $V$-reduction is a reduced commutative ring. In this case, the proof of equivalence $C\simeq \W(C/V)$ can be understood using an observation of \cite{BLM}. The map $C\rightarrow C/V$ lifts to a map of $\delta$-algebras $f: C\rightarrow \W(C/V)$ by the universal property of Witt vectors as a cofree $\delta$-ring. By $p$-torsion freeness and the fact that the ring $C/V$ is reduced, we see that this map is also compatible with $V$. Indeed, we have the equation:

$$
F(f(V(x))) = f(FV(x))=f(px)=pf(x) =FV(f(x))
$$
in $\W(C/V)$. And again, since the ring $C/V$ is reduced, the Frobenius $F: \W(C/V) \rightarrow \W(C/V)$ is injective, therefore $$f(V(x)) = V(f(x))$$ for any $x\in \W(C/V)$. Modulo $V$, the map $f$ is an equivalence, it follows by induction from $V$-completeness that it is an equivalence.

Notice the difficulty with this argument in the derived setting. For any derived $\delta$-Cartier ring $C$, the map $C\rightarrow C/V$ gives a map of $\delta$-algebras $C\rightarrow \W(C/V)$, however, a-priori one does not know that this is a map of \textbf{derived $\delta$-Cartier rings}. 
\end{rem}

\begin{ex}

Let us use the equivalence of Theorem \ref{Witt_equivalence_1} to give another computation of the free $\delta$-Cartier ring on one generator $\mathbb{Z}\{x\}[[V]]$. By the equivalence of Theorem \ref{Witt_equivalence_1}, we have:

$$
 \mathbb{Z}\{x\}[[V]]  \simeq \W(   \mathbb{Z}\{x\}[[V]] /V   ) \simeq  \W( \mathbb{Z}\{x\}  ) \simeq $$

$$ \simeq \W(\mathbb{Z}[x_{1},x_{2},...]) \simeq \W(\otimes_{i=1}^{\infty} \mathbb{Z}[x_{i}]) \simeq \widehat{\boxtimes}_{i=1}^{\infty} \W(\mathbb{Z}[x_{i}])  )
$$
Consider  the polynomial algebra $\mathbb{Z}[x]$ endowed with a lift of Frobenius given by $F(x)=x^{p}$. It is an immediate consequence of admitting a lift of Frobenius that $\W(\mathbb{Z}[x]) \simeq (\mathbb{Z}[x],F(x)=x^{p})[[V]]$. Therefore, in the end we get the following formula:

$$
 \mathbb{Z}\{x\}[[V]]  \simeq \widehat{\boxtimes}_{i=1}^{\infty} (\mathbb{Z}[x_{i}], F(x_{i}) = x_{i}^{p}) [[V]].
$$
Notice that the equivalence

$$
(\mathbb{Z}[x_{1},x_{2},...],F(x_{i})=x_{i}^{p}+px_{i+1})[[V]] \simeq  (\mathbb{Z}[x_{1},x_{2},...], F(x_{i})=x_{i}^{p} )[[V]],
$$
 is consistent with equivalence of Theorem \ref{Witt_equivalence_1} as both lifts of Frobenii coincide modulo $p$. 
\end{ex}

\subsection{Application: $p$-complete perfect derived $\delta$-rings.}

Theorem \ref{Witt_equivalence_1} recovers some classical statements and their derived analogues. An example of this is an equivalence between perfect (derived) $\mathbb{F}_{p}$-algebras and perfect (derived) $p$-complete $\delta$-rings. In the discrete case, the material below is quite classical, and we refer the reader, for instance, to Serre's book \cite[Chapter II, §5,Theorem 5]{Ser}. This was generalized to the derived setting by B.Antieau in \cite[Proposition 4.22(a)]{A23}.

\begin{defn}

Below we define the notion of perfectness in three different setting.

\begin{enumerate}
\item A derived $\mathbb{F}_{p}$-algebra $A$ is \textbf{perfect} is the derived Frobenius map $\varphi: A \rightarrow A$ is an equivalence. We let $\DAlg_{\mathbb{F}_{p}}^{\perf}$ be the $\infty$-category of perfect derived $\mathbb{F}_{p}$-algebras. This is a full subcategory of $\DAlg_{\mathbb{F}_{p}}$.

\item A derived $\delta$-ring $A\in \delta\text{-}\DAlg_{\mathbb{Z}}$ is \textbf{perfect} is the lift of Frobenius $F: A \rightarrow A$ is an equivalence. We denote $\delta\text{-}\DAlg^{\perf}_{\mathbb{Z}} \subset \delta\text{-}\DAlg_{\mathbb{Z}}$ be the full subcategory spanned by perfect derived $\delta$-algebras. The inclusion is a localization. The left adjoint to the inclusion is given by the formula

\begin{equation}\label{colimit_perfection}
A_{\perf}:=\colim ( A \xymatrix{ \ar[r]^-{\varphi}&} A  \xymatrix{ \ar[r]^-{\varphi}&}  A   \xymatrix{ \ar[r]^-{\varphi}&}   ...  )
\end{equation}
for any derived $\delta$-algebra $(A,\varphi)$.

\item We define the $\infty$-category of \textbf{perfect $V$-complete derived $\delta$-Cartier rings} $\delta\text{-}\widehat{\Cart}\DAlg_{\mathbb{Z}_{p}}^{\perf}$ as the full subcategory of $ \delta\text{-}\widehat{\Cart}\DAlg_{\mathbb{Z}_{p}}$ corresponding to $\DAlg_{\mathbb{F}_{p}}^{\perf} \subset \DAlg_{\mathbb{F}_{p}}$ under the equivalence $\W: \DAlg_{\mathbb{F}_{p}} \simeq \delta\text{-}\widehat{\Cart}\DAlg_{\mathbb{Z}_{p}}$. One consequence of the fact that the composition $\W: \DAlg_{\mathbb{F}_{p}} \simeq \delta\text{-}\widehat{\Cart}\DAlg_{\mathbb{Z}_{p}} \rightarrow \delta\text{-}\DAlg_{\mathbb{Z}_{p}}$ is equivalent to the classical Witt vectors functor, is that for any $A\in \DAlg_{\mathbb{F}_{p}}^{\perf}$, the Witt vectors Frobenius $F: \W(A) \rightarrow \W(A)$ is an equivalence. In fact, the full subcategory $\delta\text{-}\widehat{\Cart}\DAlg_{\mathbb{Z}_{p}}^{\perf} \subset \delta\text{-}\widehat{\Cart}\DAlg_{\mathbb{Z}_{p}}$ consists of objects whose underlying derived $\delta$-algebra is perfect.
\end{enumerate}

\end{defn}

\begin{rem}
Note that since we are working over $\mathbb{Z}_{p}$, for any $A\in  \delta\text{-}\widehat{\Cart}\DAlg_{\mathbb{Z}_{p}}$, the Frobenius $F: A \rightarrow A$ is an endomorphism of $\delta$-Cartier rings (it commutes with the $\delta$-structure and the Verschiebung). Consequently, we can characterize the inclusion $\delta\text{-}\widehat{\Cart}\DAlg_{\mathbb{Z}_{p}}^{\perf} \subset \delta\text{-}\widehat{\Cart}\DAlg_{\mathbb{Z}_{p}}$ as the localization at the class of morphisms $F: A \rightarrow A$ provided by the Frobenius map for any $A\in \delta\text{-}\widehat{\Cart}\DAlg_{\mathbb{Z}_{p}}$. In particular, the left adjoint of the inclusion $\delta\text{-}\widehat{\Cart}\DAlg_{\mathbb{Z}_{p}}^{\perf} \subset \delta\text{-}\widehat{\Cart}\DAlg_{\mathbb{Z}_{p}}$  is given by the composition of the same formula as \ref{colimit_perfection}, followed by $p$-completion.
\end{rem}

We now prove the following Lemma. 

\begin{Lemma}\label{delta_perfection}
The forgetful functor $\delta\text{-}\widehat{\Cart}\DAlg_{\mathbb{Z}_{p}}^{\perf} \rightarrow \delta\text{-}\DAlg_{\mathbb{Z}_{p}}^{\perf,\wedge} $ is an equivalence.
\end{Lemma}

\begin{proof}
Indeed, unwinding the definitions, the left adjoint of this functor is given by sending $A \in \delta\text{-}\DAlg_{\mathbb{Z}_{p}}^{\perf,\wedge}$ to $A[V]^{\wedge}_{\perf}$ (the argument below shows that $A[V]_{\perf}$ is already $p$-complete, and hence $V$-complete too). Computing $A[V]_{\perf}$ by the colimit perfection formula, we get that $$A[V]_{\perf} \simeq \underset{F}{\colim} (A[V] \rightarrow A[V] \rightarrow ...) \simeq \underset{F}{\colim} (\underset{n\geq 0}{\bigoplus} V^{n}A \rightarrow \underset{n\geq 0}{\bigoplus}V^{n}A \rightarrow ... ) \simeq $$

$$\simeq \underset{n}{\colim} \:\underset{F}{\colim} (\bigoplus_{i=0}^{n}V^{i}A \rightarrow \bigoplus_{i=0}^{n}V^{i}A \rightarrow ...) $$
Since for any direct summand $V^{i}A \subset \oplus_{j=0}^{n} V^{j}A$, the $i$-th power of Frobenius lands in $A$, it follows that the colimit above is equivalent to $A$, and hence $A$ has the structure of a perfect derived $\delta$-Cartier ring with $V=p\varphi^{-1}$. Consequently, the functor $(-)[V]_{\perf}$ is the inverse of the forgetful functor.
\end{proof}

Applying Lemma \ref{delta_perfection}, we obtain the next Theorem as a consequence of Theorem \ref{Witt_equivalence_relative}.

\begin{Theor}
There is an equivalence of $\infty$-categories.

$$
\xymatrix{(-)/p : \delta\text{-}\DAlg_{\mathbb{Z}_{p}}^{\perf,\wedge}  \ar@<+1.0ex>[rrr]_-{\sim} &&& \ar@<+1.0ex>[lll] \DAlg_{\mathbb{F}_{p}}^{\perf} : \W, }
$$
\end{Theor}

\end{document}

%% file: lang_en.tex
\usepackage[T1]{fontenc}
\usepackage[utf8]{inputenc}
\usepackage[english]{babel}
\usepackage[autostyle, english = american]{csquotes}
\MakeOuterQuote{"}

\usepackage{amsthm}
\theoremstyle{definition}
\newtheorem{Theor}{Theorem}[subsection]
\newtheorem*{thm*}{Theorem}
\newtheorem{Lemma}[Theor]{Lemma}

\newtheorem*{prop*}{Proposition}

\newtheorem{warning}[Theor]{Warning}
\newtheorem{construction}[Theor]{Construction}

\newtheorem{variant}[Theor]{Variant}

\newtheorem{defn}[Theor]{Definition}
\newtheorem{lem}[Theor]{Lemma}
\newtheorem{cor}[Theor]{Corollary}
\newtheorem{prop}[Theor]{Proposition}
\newtheorem{ex}[Theor]{Example}
\newtheorem{rem}[Theor]{Remark}
\newtheorem{question}[Theor]{Question}

\usepackage{color}

\definecolor{note_color}{rgb}{0.0,0.7,0.0}

%% file: hat.tex
\usepackage{mathrsfs}
\usepackage{latexsym}
\usepackage{graphicx}
\usepackage{amscd,amssymb,amsmath,amsbsy,amsfonts}
\usepackage{amsmath}
\usepackage{mathabx}
\usepackage[mathscr]{eucal}
\usepackage[all]{xy}
\usepackage{color}
\definecolor{reference}{rgb}{0.20,0.36,0.74}
\definecolor{citation}{rgb}{0,.40,.80}
\usepackage[colorlinks,linktocpage,urlcolor=blue,linkcolor=reference,citecolor=citation]{hyperref}
\usepackage{verbatim}
\usepackage{enumerate}
\usepackage{cleveref}
\crefformat{equation}{(#2#1#3)}
\usepackage{tikz}
\usepackage{graphicx}
\usepackage[left=2cm, top=2cm, right=2cm, bottom=2cm]{geometry}

\DeclareFontFamily{OT1}{pzc}{}
\DeclareFontShape{OT1}{pzc}{m}{it}{<-> s * [1.200] pzcmi7t}{}
\DeclareMathAlphabet{\mathpzc}{OT1}{pzc}{m}{it}

\usepackage{enumitem}

\newlength{\bibitemsep}
\newlength{\bibparskip}
\let\oldthebibliography\thebibliography
\renewcommand\thebibliography[1]{%
  \oldthebibliography{#1}%
  \setlength{\parskip}{\bibitemsep}%
  \setlength{\itemsep}{\bibparskip}%
}

\makeatletter
\DeclareRobustCommand{\emdef}{%
  \@nomath\em \if b\expandafter\@car\f@series\@nil
  \normalfont \else \bfseries \fi}
\makeatother

\newcommand{\fcat}{\mathscr}

\newcommand{\Set}{ {\fcat{S}\mathrm{et}} }

\renewcommand{\phi}{\varphi}
\renewcommand{\epsilon}{\varepsilon}

\DeclareMathOperator{\Sym}{Sym}
\newcommand{\Mod}{\mathrm{Mod}}

\DeclareMathOperator{\gr}{gr}
\DeclareMathOperator{\rk}{rk}

\newcommand{\perf}{{\mathrm{perf}}}

\DeclareMathOperator{\Frob}{Frob}

\DeclareMathOperator{\Hom}{Hom}

\DeclareMathOperator{\Tor}{Tor}

\DeclareMathOperator{\End}{End}

\DeclareMathOperator{\Map}{Map}
\DeclareMathOperator{\Fun}{Fun}

\DeclareMathOperator*{\colim}{colim}
\DeclareMathOperator*{\fib}{fib}

\DeclareMathOperator*{\cofib}{cofib}
\DeclareMathOperator*{\coker}{coker}

\DeclareMathOperator{\Alg}{Alg}
\DeclareMathOperator{\CAlg}{CAlg}

\DeclareMathOperator{\Tot}{Tot}

\DeclareMathOperator{\Env}{Env}

\mathchardef\mdef="2D